\documentstyle{amsppt}
\TagsOnRight
\loadbold
\mathsurround=1.2pt

\define\itm#1#2{\smallskip\noindent{\rlap{\rom{(#1)}}\qquad}#2\smallskip}
\define\Y#1{{\kern1pt{:}\kern0.5pt#1\kern-0.3pt{:}\kern0.4pt}}
\define\sbh#1{\smallskip\noindent#1\quad}


\overfullrule=0pt
\TagsOnRight

\def\myhookr#1{\mathrel{%
	     \mathrel{\raise .3ex\hbox{$\subset$}}%
	     \joinrel\mathrel{\vrule height.35ex depth-.25ex width#1}%
	     \joinrel\mathrel{\lower .3ex\hbox{$\longrightarrow$}}}}

\def\myhookl#1{\mathrel{\lower.3ex\hbox{$\longleftarrow$}}%
		\joinrel\mathrel{\vrule height.35ex depth-.25ex width#1}%
		\joinrel\mathrel{\raise.3ex\hbox{$\supset$}}}
\def\sqr#1#2{{\vcenter{\vbox{\hrule height.#2pt
        \hbox{\vrule width.#2pt height#1pt \kern#1pt
          \vrule width.#2pt}
        \hrule height.#2pt}}}}
\def\square{\mathchoice\sqr34\sqr34\sqr{2.1}3\sqr{1.5}3}

\define\dd{{\buildrel{_{\leftrightarrow}}\over\Delta}^{\!}}
\define\Per{{\Cal P}\text{\it er}}
\define\K{{\bold K}}
\define\DD{{\bold D}}
\define\Dt{{\bold {D^t}}}
\define\G{\Gamma^\diamondsuit}
\define\I{{I}}
\define\J{{\Cal J}}
\define\PP{{P}}
\define\la{{|\lambda|}}
\define\ma{{|\mu|}}

\define\iso{{\,{\buildrel{{\sim}}\over\longrightarrow}\,}}

\define\Ug{{{\text{\rm U}}\!{\frak g}}}
\define\Uh{{{\text{\rm U}}\!{\frak h}}}
\define \Spec {\operatorname{Specm}} 
\define \Res {\operatorname{Res}}
\define \id {\operatorname{Id}}
\define \ad {\operatorname{ad}}
\def\g{{\frak g}}
\def\D{{\Bbb D}}
\def\C{{\Bbb C}}
\define \Mod {\text{Mod}}
\define \HC {\text{HC}}
\define \rk {\operatorname{rk}}
\define \reg {{\text{\rm reg}}}  
\define \pr {\operatorname{pr}}
\define \Hom {\operatorname{Hom}}
\define \Hol {\operatorname{\Cal Hol}}
\define \End {\operatorname{End}}

\baselineskip = 18pt
\magnification=1100

\topmatter

\title 
WALL-CROSSING FUNCTORS AND ${\Cal D}$-MODULES
\endtitle
\author Alexander Beilinson and Victor Ginzburg
\endauthor

\abstract
We study Translation functors  and  Wall-Crossing functors on
infinite dimensional representations of  a complex semisimple
Lie algebra using ${\Cal D}$-modules. This functorial  machinery is then
used to prove the Endomorphism-theorem and the Structure-theorem,
two important results  established earlier by  W. Soergel in a 
totally different way. Other applications to the category ${\Cal O}$ of
Bernstein-Gelfand-Gelfand are given, and some conjectural
relationships between Koszul duality, Verdier duality and
convolution functors are discussed. A geometric interpretation of
tilting modules is given.
\endabstract

\address A.B.: Department of Mathematics, University of Chicago, Chicago IL 
60637, USA\newline
$\hphantom{x}\qquad\qquad{}^{\text{E-mail: {\bf
sasha\@math.uchicago.edu}}}$\endaddress
\address V.G.: Department of Mathematics, University of Chicago, Chicago IL 
60637, USA\newline
$\hphantom{x}\qquad\qquad{}^{\text{E-mail: {\bf
ginzburg\@math.uchicago.edu}}}$
\endaddress
\endtopmatter


\document

\vskip 5mm

\centerline{\bf Table of Contents}
\vskip 10mm

\qquad${}^{\text{\bf 1.{ $\;$} Translation functors for the extended enveloping
algebra}}_{\text{\bf 2.{ $\;$} Translation functors via
${\Cal D}$-modules}}$

\vskip 1pt

\qquad${}^{\text{\bf 3.{ $\;$} Properties of the Translation
functors}}_{\text{\bf 4.{ $\;$} Applications to the category ${\Cal
O}$}}$

\vskip 1pt

\qquad${}^{\text{\bf 5.{ $\;$} Convolution of Harish-Chandra 
modules}}_{\text{\bf 6.{ $\;$} Projective 
functors ${\vphantom{X}}^{}$}}$

\vskip 6mm

\subheading{0.\ Introduction}
\vskip .15in

We will be concerned here with infinite dimensional
representations  of a complex semisimple Lie algebra ${\frak g}$.

In more detail, let $\Ug$ be the universal enveloping algebra of ${\frak
g}$, and  let $Z({\frak g})$ be the center of $\Ug$. We consider the
category of $\Ug$-modules annihilated by a
great enough (unspecified) power of a maximal ideal in $Z({\frak
g})$. It has been observed
by many people during the 70's, see e.g., [W], [LW], [J], and [Z],
that various results can be
usefully  transferred between the categories  
corresponding to two different maximal
ideals,
 using tensor
products with finite-dimensional representations.
The most relevant for us is the work of Jantzen [J], who
 introduced certain functors between the two
categories, called {\it translation functors}. Jantzen showed 
that if both  maximal
ideals satisfy certain regularity and integrality conditions,
then the translation functor establishes an equivalence of the
two categories. If one of the two ideals is regular while the other is
not, the corresponding translation 
functor is no longer an equivalence. The composition
of the translation functor that
sends the category at a regular maximal ideal
to the category at a non-regular maximal ideal with the translation
functor acting in the opposite direction is called a {\it wall-crossing}
functor. The terminology stems from the identification of (integral)
maximal ideals of $Z({\frak g})$
with Weyl group orbits
 in the weight lattice of the maximal torus. 
 Non-regular ideals correspond to the orbits contained in
 the union of walls
of the Weyl chambers.

In this paper we will be mainly interested in the "most singular" case,
where the non-regular maximal ideal corresponds to the fixed point
of the Weyl group, that is the point contained in all the walls.
Our study was partly motivated  by trying to understand two
important results, the 
"Endomorphism-theorem" and the "Structure-theorem",
 proved
by W. Soergel [S1] in the course of the proof of the Koszul duality
conjecture, see [BG], [BGS]. Soergel's argument was very clever,
 but rather technical. A shorter
proof of the Endomorphism-theorem was found by J. Bernstein [Be].
An alternative entirely geometric approach to the
Structure-theorem, based on perverse
sheaves and leaving almost all of the Representation theory aside
was given in [Gi2]. We propose below  new
proofs of both theorems in the framework of Representation theory,
and based on the technique of ${\Cal D}$-modules. An advantage of
our approach is that the bulk of the argument goes through in a quite
general setting, while the proofs in [S1] and [Be] relied heavily on
some special
features of the category ${\Cal O}$ right from the beginning.
One of our goals is to convince the reader (or at least ourselves)
that the results in question are not so "hard", and that they
follow quite naturally from the ${\frak g}$-{\it module} $\longleftrightarrow$
${\Cal D}$-{\it module} correspondence combined with some basic functoriality
properties.

It would be  very interesting to apply  the general results
obtained in section 3 below to the category of Harish-Chandra modules instead
of the category ${\Cal O}$. It is also tempting
to extend our results to representations of affine Lie algebras
(cf. [FrMa]). We remark that the ${\Cal D}$-module
approach of the present paper
provides natural and simple proofs of most of the
results concerning {\it tilting} modules (over the 
{\it finite-dimensional} semisimple Lie algebra $\g$),
  proved by Soergel in [S2] 
for tilting modules over an {\it affine} Lie algebra.  The main
problem of extending our results to the affine case is that the 
category at the regular  {\it positive} level is known [KT1] to
be related to ${\Cal D}$-modules on one affine flag manifold
(which is the union of cells of finite {\it co}dimension);
the 
category at the regular  {\it negative} level is known [KT2] to
be related to ${\Cal D}$-modules on another affine flag manifold
(which is the union of cells of finite dimension), and
the 
category at the "most singular" = {\it critical} level is expected
[FM] to
be related to ${\Cal D}$-modules on yet another,
so-called {\it periodic} flag manifold,  
(which is the union of cells of {\it semi-infinite} dimension).
It is therefore especially intriguing how to find a
geometric construction of
translation functors to and from the "most singular"
 maximal ideal.

We now briefly outline the contents of the paper. In \S 1
we introduce a Galois extension of $\Ug$ with the Galois group $W$,
the Weyl group. We compare modules over $\Ug$ with those over
the Galois extension. Section 2 provides a self-contained
exposition of translation functors from the geometric view-point.
Although some of the results of this section are undoubtedly known
to experts, we could not find a relevant reference in the literature.
Section 3 is the heart of the paper containing our main results
about wall-crossing functors. 
In \S 4 we 
prove
the {\it Endomorphism-theorem} and the {\it Structure-theorem} of Soergel,
give some applications to the category
${\Cal O}$ of Bernstein-Gelfand-Gelfand [BGG].

Section 5 is devoted to various
convolution functors. A  convolution functor was first introduced in the
geometric setting of perverse sheaves back in 1980, independently
by Beilinson-Bernstein, Brylinski, Lusztig, MacPherson and others,
in the course of the proof of the Kazhdan-Lusztig conjecture.
An algebraic counterpart of that convolution in terms of Harish-Chandra
bi-modules was considered in [Gi]. It was claimed in [Gi] (without
proof) that the algebraic convolution of Harish-Chandra
bi-modules goes into the geometric convolution under the localization
functor of [BB]. This claim turns out to be not quite correct, and in $\S 5$
we establish
a precise  relationship between the algebraic and the
geometric convolution functors, respectively. We also propose
conjectural  relationships (see Conjecture 5.18 and Theorem 5.24)
between the various convolutions,
Koszul duality of [BGS], and the Verdier duality functors.

In $\S 6$ we express projective functors defined in [BeGe]
in terms of the convolution functors studied in the previous section.
We establish a direct connection
between projective functors and projective Harish-Chandra
bi-modules (somewhat analogous to a result of [BeGe]
that motivated the name "projective functor").
We discuss some applications. In particular,
we introduce {\it tilting} Harish-Chandra
modules, establish their relation to the tilting
objects of the categery ${\Cal O}$, and derive
their basic properties using convolution functors.

The present paper grew out of
 an unpublished chapter of one of the (many) preliminary
versions of [BGS], written in 1992. 
That chapter was not directly related to the subject
of [BGS], and we decided to publish it separately to keep the size
of [BGS] to a minimum.
\vskip .2in

\subheading{1.\ Translation functors for the extended enveloping algebra}
\vskip .15in

Throughout the paper we fix ${\frak g}$, a complex semisimple Lie
algebra. We write ${\frak h}$ for the corresponding
{\it abstract} Cartan subalgebra. The reader should be warned that
${\frak h}$ is not a subalgebra of ${\frak g}$; it is defined, see
e.g. [CG, p.137],
 as the
quotient of a Borel subalgebra modulo its nil-radical, and this quotient
is independent
of any choices.
Recall further, see {\it loc. cit.}, that 
${\frak h}^*$, the dual space, comes equipped with a root system
$R\subset {\frak h}^*$ which has a preferred choice of simple roots.
Thus there is a well-defined element $\rho\in{\frak h}^*$= {\it half-sum
of positive roots}. An element $\lambda\in {\frak h}^*$ is called
$\rho$-{\it dominant} if $\langle\lambda+\rho,\check{\alpha}\rangle
\ge 0\,,$ for every positive coroot $\check{\alpha}$.

Let
$W$ be the Weyl group of ${\frak g}$, the Coxeter group generated by reflections
with respect to roots. The group $W$ acts naturally on
${\frak h}$ and on ${\frak h}^*$. In this paper we will always use
the so-called
{\it dot-action} of $W$ on ${\frak h}^*$.  
The dot-action of $w\in W$
is obtained by "twisting"
the standard $w$-action:  ${\frak h}^*\ni\lambda\mapsto w(\lambda)$
as follows: $\lambda\mapsto w\cdot\lambda :=w(\lambda+\rho)-\rho\,.$
Thus, the point $(-\rho)$ is the unique $W$-fixed point of
the dot-action. 
Given $\lambda \in {\frak h}^*$, let 
$W_\lambda\subset W$ denote the isotropy group, and
$\la \simeq W/W_\lambda$  
the $W$-orbit of $\lambda$ with respect 
to the dot-action. 

The dot-action on ${\frak h}^*$ induces a $W$-action on the polynomial
algebra ${\Bbb C}[{\frak h}^*]$.
Let  ${\Bbb C}[{\frak h}^*]^W \subset {\Bbb C}[{\frak h}^*]$  be 
the subalgebra of $W$-invariant polynomials on  ${\frak h}^*$. We have 
$\Spec {\Bbb C}[{\frak h}^*]^W = {\frak h}^*/W$,  where  $\Spec$ stands for the maximal 
spectrum of a commutative  ${\Bbb C}$-algebra.
 Given 
$\lambda \in {\frak h}^*$, we let  $\J_\lambda \subset 
{\Bbb C}[{\frak h}^*]$  denote the maximal ideal of all polynomials vanishing 
at $\lambda$. We will often identify ${\Bbb C}[{\frak h}^*]$ with $S{\frak h}$, the
Symmetric algebra on ${\frak h}$.

Let  $Z({\frak g})$  be the center of the universal
enveloping algebra   $\Ug$.  
The Harish-Chandra isomorphism $Z({\frak g})
\simeq S{\frak h}^W$
gives rise to the following composition of algebra morphisms:
$$\chi:
Z({\frak g}) \iso S{\frak h}^W \hookrightarrow S{\frak h} \iso
 {\Bbb C}[{\frak h}^*]\,.$$
We will identify  $Z({\frak g})$  with the image of $\chi$,
 the subalgebra ${\Bbb C}[{\frak h}^*]^W$ of dot-invariant polynomials. Thus, the 
map  $\Spec {\Bbb C}[{\frak h}^*] \to \Spec Z({\frak g})$  induced by 
the imbedding $\chi$ can (and will)
be identified with the projection:   $\pi: {\frak h}^* \twoheadrightarrow 
{\frak h}^*/W$.  
Given $\lambda \in {\frak h}^*$, we will often view the 
$W$-orbit $\la$ as a point of the orbi-space
${\frak h}^*/W$, and we let  $\I_{\la} =\pi^{-1}(\J_\lambda)
\in \Spec Z({\frak g})$
denote the corresponding maximal ideal in $Z({\frak g})$. 

\vskip 3pt 
\proclaim{Definition}
The extended enveloping algebra is defined by  $\widetilde{U} := \Ug 
\otimes_{Z({\frak g})} {\Bbb C}[{\frak h}^*]$. 
\endproclaim
\noindent
 Thus  $\widetilde{U}$  contains 
both $\Ug$ and ${\Bbb C}[{\frak h}^*]$  as subalgebras and the subalgebra 
${\Bbb C}[{\frak h}^*]$  coincides with the center of  $\widetilde{U}$.  
Furthermore, the $W$-action on ${\Bbb C}[{\frak h}^*]$  gives rise to a 
$W$-action on $\widetilde{U}$  such that  $\widetilde{U}^W = \Ug$.  

Given  $\lambda \in {\frak h}^*$, let  $\Mod_{\la}(\Ug)$, resp. 
$\Mod_\lambda(\widetilde{U})$, denote the category of finitely-generated 
$\Ug$-modules (resp. $\widetilde{U}$-modules) $M$ such that 
$\I^n_{\la} \cdot M= 0$, resp. $\J^n_\lambda \cdot M = 0$,
for great enough  $n = n(M) \gg 0$.  Restricting $\widetilde{U}$-modules to 
$\Ug$-modules yields an exact functor  $\Res_\lambda: \Mod_\lambda(\widetilde{U}) 
\longrightarrow \Mod_{\la}(\Ug)$.  

Fix $\lambda \in {\frak h}^*$. We introduce 
the intermediate algebra of  $W_\lambda$-invariants:
${\Bbb C}[{\frak h}^*]^W \subset{\Bbb C}[{\frak h}^*]^{W_\lambda} \subset 
{\Bbb C}[{\frak h}^*]$. We have
$$\aligned
\widetilde{U}\; \simeq\;
\Ug\otimes_{{{\Bbb C}[{\frak h}^*]^W}}{\Bbb C}[{\frak h}^*]&
\;\simeq\;
\Ug\otimes_{{{\Bbb C}[{\frak h}^*]^{W}}}
{\Bbb C}[{\frak h}^*]^{W_\lambda}\otimes_{{{\Bbb C}[{\frak h}^*]^{W_\lambda}}}
{\Bbb C}[{\frak h}^*]\\
&\qquad \qquad\simeq\;\widetilde{U}^{W_\lambda}
\otimes_{{{\Bbb C}[{\frak h}^*]^{W_\lambda}}}
{\Bbb C}[{\frak h}^*]
\endaligned
\tag 1.1
$$

Geometrically, introducing the intermediate algebra
${\Bbb C}[{\frak h}^*]^W \hookrightarrow
{\Bbb C}[{\frak h}^*]^{W_\lambda} \hookrightarrow
{\Bbb C}[{\frak h}^*]$ corresponds to the factorisation of the projection
${\frak h}^* \twoheadrightarrow {\frak h}^*/W$ as the composition
${\frak h}^* \twoheadrightarrow {\frak h}^*/W_\lambda
\twoheadrightarrow {\frak h}^*/W$. Let
 $\J_\lambda^{W_\lambda} := 
{\Bbb C}[{\frak h}^*]^{W_\lambda} \cap \J_\lambda$
be  the corresponding maximal ideal.
Observe that the projection
${\frak h}^*/W_\lambda \twoheadrightarrow
 {\frak h}^*/W$  is unramified over $\la$.
Hence, the  $\I_{\la}$-adic completion of the algebra 
$Z({\frak g})$  is isomorphic {\it canonically} to the  $\J_\lambda^{W_\lambda}$-adic
completion of the algebra  ${\Bbb C}[{\frak h}^*]^{W_\lambda}$.  This yields a 
canonical algebra isomorphism between the  $\I_{\la}$-adic 
completion of  $\Ug$  and the  $\J_\lambda^{W_\lambda}$-adic completion of 
$\widetilde{U}^{W_\lambda}$.  Now, let  $\Mod_\lambda(\widetilde{U}^{W_\lambda})$
be the category of finitely-generated $\widetilde{U}^{W_\lambda}$-modules 
annihilated by a great enough power of the ideal  $\J_\lambda^{W_\lambda}$. For 
any $M \in \Mod_\lambda(\widetilde{U}^{W_\lambda})$, the 
$\widetilde{U}^{W_\lambda}$-action on $M$ can be uniquely extended, by continuity,
to an action of the  $\J_\lambda^{W_\lambda}$-adic completion; similarly, the 
$\Ug$-action on any  $M \in \Mod_{\la}(\Ug)$  can be 
extended to an action of the  $\I_{\la}$-adic completion.
The completions being isomorphic, we obtain

\proclaim{Lemma 1.2}  The $\Ug$-action on any module  
$M \in \Mod_{\la}(\Ug)$  can be extended canonically to an 
$\widetilde{U}^{W_\lambda}$-action so that the following
restriction functor is an equivalence of categories:
$$
\Mod_\lambda(\widetilde{U}^{W_\lambda}) \longrightarrow
\Mod_{\la}(\Ug)\;.
\quad\square
$$
\endproclaim

\noindent
We will often view  $\Ug$-modules as  $\widetilde{U}^{W_\lambda}$-modules
via the Lemma.

\medpagebreak

Observe next that for the projection  $\pi: {\frak h}^* \to {\frak h}^*/W_\lambda$
we have  $\pi^{-1} \circ \pi(\lambda) = \lambda$.  It follows that for any 
$M \in \Mod_\lambda(\widetilde{U}^{W_\lambda})$  one has:  
$\widetilde{U} \otimes_{\widetilde{U}^{W_\lambda}}M \in \Mod_\lambda(\widetilde{U})$.
The functor  $\widetilde{U} \otimes_{\widetilde{U}^{W_\lambda}}(\cdot )$ is 
clearly the left adjoint of the restriction functor  $\Mod_\lambda(\widetilde{U}) 
\longrightarrow  \Mod_\lambda(\widetilde{U}^{W_\lambda})$.  Thus, decomposing 
the functor  $\Res_\lambda: \Mod_\lambda(\widetilde{U})\longrightarrow 
\Mod_{\la}(\Ug)
$  as the composition $\Mod_\lambda(\widetilde{U}) 
\to \Mod_\lambda(\widetilde{U}^{W_\lambda}) \to
\Mod_{\la}(\Ug)$  and taking adjoints, we get 

\proclaim{Lemma 1.3}  The functor $M \mapsto \widetilde{U} 
\otimes_{\widetilde{U}^{W_\lambda}}M$  is the left adjoint of the functor
$\Res_\lambda$, where  $M$  is viewed as a $\widetilde{U}^{W_\lambda}$-module
via Lemma 1.2.\quad $\square$
\endproclaim

\medpagebreak

\medpagebreak

We will call a weight $\lambda \in {\frak h}^*$

\noindent
$\bullet\quad${\it integral} if
$\langle\lambda,\check{\alpha}\rangle\in{\Bbb Z}\,,$
for any coroot $\check{\alpha}$;

\noindent
$\bullet\quad${\it regular} if  $W_\lambda = 
\{ 1\}$.   
\vskip 3pt

 Let  $E$  be a 
finite-dimensional ${\frak g}$-module and $M \in \Mod_{\la}(\Ug)$,
$\lambda \in {\frak h}^*$.  Then, the $\Ug$-module  $E \otimes_{\Bbb C}M$
is annihilated by an ideal in $Z({\frak g})$  of finite codimension [Ko].  Hence
there is a canonical finite direct sum decomposition:
$$
E \otimes M = {\bigoplus_{\ma\in {\frak h}^*/W}}\text{pr}_\ma(E \otimes 
M)\,,\quad
\text{where}\quad \text{pr}_\ma(E \otimes M) \in \Mod_{\ma}(\Ug).
$$
Fix integral $\rho$-dominant weights  $\lambda, \mu 
\in {\frak h}^*$, and let
 $E_{\lambda - \mu}$  be an irreducible finite-dimensional 
${\frak g}$-module with extreme weight  $\lambda - \mu$.  Following [Ja], 
define the translation functor 
$$
\theta^\lambda_\mu: \Mod_{\ma}(\Ug) \longrightarrow 
\Mod_{\la}(\Ug)\quad \text{by}\quad 
\theta^\lambda_\mu M = \text{pr}_\la(E_{\lambda-\mu} \otimes M)\,.
$$
The functor  $\theta^\mu_\lambda$  is both the left and the right adjoint of
$\theta^\lambda_\mu$.  The functors $\theta^\lambda_\mu$ and
$\theta^\mu_\lambda$ are both exact.

By Lemma 1.2, the functor $\theta^\lambda_\mu$  can be viewed as a functor
$\Mod_\mu(\widetilde{U}^{W_\mu}) \to \Mod_\lambda(\widetilde{U}^{W_\lambda})$.
The action of any  $z \in {\Bbb C}[{\frak h}^*]^{W_\mu}$  on  
$M \in \Mod_\mu(\widetilde{U}^{W_\mu})$ is a $\Ug$-module 
endomorphism,  hence induces by functoriality an endomorphism 
$\theta^\lambda_\mu(z): \theta^\lambda_\mu M \to \theta^\lambda_\mu M$.
One can
 describe this endomorphism in terms of the algebra automorphism  
$T_{\lambda -\mu}: {\Bbb C}[{\frak h}^*] \to {\Bbb C}[{\frak h}^*]$  induced by 
the affine translation:  $x \longmapsto x + (\lambda - \mu),\ \  x \in
{\frak h}^*$,
as follows.

Assume in addition that  $W_\lambda \subset W_\mu$.  Then, the automorphisms 
$T_{\pm(\lambda -\mu)}$  preserve the subalgebra  ${\Bbb C}[{\frak h}^*]^{W_\lambda}$. 
We also have ${\Bbb C}[{\frak h}^*]^{W_\mu} \subset {\Bbb C}[{\frak h}^*]^{W_\lambda}$, 
hence  $T_{\pm(\lambda -\mu)}{\Bbb C}[{\frak h}^*]^{W_\mu}$
$\subset 
{\Bbb C}[{\frak h}^*]^{W_\lambda}$.  In the next section we will prove
the following result of Soergel [S1,~Thm. 8]
using the  ${\Cal D}$-module approach.

\proclaim{Propositon 1.4}  Let  $z \in {\Bbb C}[{\frak h}^*]^{W_\mu}$. 
Then:  

\item{(i)}  For any  $M \in \Mod_\mu(\widetilde{U}^{W_\mu})$  and  
$m \in \theta^\lambda_\mu M$  we have:  $z \cdot m = 
[T_{\lambda -\mu} \theta^\lambda_\mu(z)] \cdot m$;

\item{(ii)}  For any  $M \in \Mod_\lambda(\widetilde{U}^{W_\lambda})$  and  
$m \in \theta^\mu_\lambda M$  we have:  $z \cdot m = 
[\theta_\lambda^\mu(T_{\lambda -\mu} z)] \cdot m$.
\endproclaim

From now on,  fix integral $\rho$-dominant weights $\lambda,\mu \in {\frak h}^*$,
such that  $W_\lambda \subset W_\mu,$  and set
$\theta^+ := \theta^\lambda_\mu$, and $\theta^- := \theta^\mu_\lambda$.  We are 
going to extend the functors  $\theta^\pm$ to $\widetilde{U}$-modules.

First, consider the composition of functors  
$$\Mod_\lambda(\widetilde{U})\ 
{\buildrel {{\Res_\lambda}}\over\longrightarrow}\ \Mod_{\la}(\Ug)\ 
{\buildrel {\theta^-}\over\longrightarrow}\ \Mod_{\ma}(\Ug)\,.$$
  The
action of an element  $a \in {\Bbb C}[{\frak h}^*] \subset \widetilde{U}$  on 
$M \in \Mod_\lambda(\widetilde{U})$  induces, by functoriality, an 
endomorphism  $\theta^-(a):  \theta^-({\Res}_\lambda M) \to 
\theta^-({\Res}_\lambda M)$.  We define a ${\Bbb C}[{\frak h}^*]$-action on 
$\theta^-({\Res}_\lambda M)$  by the formula:
$$
a*m := \theta^-(T_{\mu -\lambda} a) \cdot m,\qquad m \in \theta^-({\Res}_\lambda M).
$$
When restricted to the subalgebra  ${\Bbb C}[{\frak h}^*]^{W_\mu} \subset 
{\Bbb C}[{\frak h}^*]$, this action coincides, by Proposition 1.4(ii), with the
${\Bbb C}[{\frak h}^*]^{W_\mu}$-action on $\theta^-({\Res}_\lambda M)$  
arising from the  $\widetilde{U}^{W_\mu}$-module structure.  Therefore,
combining together the 
${\Bbb C}[{\frak h}^*]$- and the  $\widetilde{U}^{W_\mu}$-actions we get
a $\widetilde{U}^{W_\mu}\otimes_{\Bbb C}{\Bbb C}[{\frak h}^*]$-action
which factors through an action of the algebra
$\widetilde{U}^{W_\mu}\otimes_{_{{\Bbb C}[{\frak h}^*]^{W_\mu}}}
{\Bbb C}[{\frak h}^*]$. In view of isomorphism (1.1) that gives
a $\widetilde{U}$-action on
$\theta^-({\Res}_\lambda M)$.  This way we obtain an exact functor 
$\widetilde{\theta}^-: \Mod_\lambda(\widetilde{U}) \to \Mod_\mu(\widetilde{U})$.

Next, consider the composition  $\Mod_\mu(\widetilde{U})\ 
{\buildrel {{\Res_\mu}}\over\longrightarrow}\ \Mod_{\ma} (\Ug)\ 
{\buildrel {\theta^+}\over\longrightarrow}\ \Mod_{\la}(\Ug)$.  For any 
$N \in \Mod_\mu(\widetilde{U})$, define a ${\Bbb C}[{\frak h}^*]$-action 
on  $\theta^+({\Res}_\mu N)$  by
$$
a*m := [T_{\lambda - \mu} \theta^+(a)] \cdot m,\quad  
m \in \theta^+({\Res}_\mu N).
$$
Again, when restricted to the subalgebra  ${\Bbb C}[{\frak h}^*]^{W_\mu}$, 
this action coincides, by Proposition 1.4(i), with the action arising from the 
natural  $\widetilde{U}^{W_\lambda}$-module structure on  $\theta^+({\Res}_\mu N)$.
This is not necessarily the case, however, for the full algebra 
${\Bbb C}[{\frak h}^*]^{W_\lambda}  \supset {\Bbb C}[{\frak h}^*]^{W_\mu}$;
the two actions might be different!  Thus, we put:
$$
\widetilde{\theta}^+_r N = \{ n \in \theta^+(\Res_\mu N)
\enspace\big\vert\quad
 a*n = a\cdot n,\quad
\forall a \in {\Bbb C}[{\frak h}^*]^{W_\lambda}\}\,,
$$
and dually,
$$
\widetilde{\theta}^+_\ell N = \theta^+(\Res_\mu N)/\{ a*n 
- a\cdot n\enspace \big\vert\quad a \in 
{\Bbb C}[{\frak h}^*]^{W_\lambda},\quad n \in \theta^+(\Res_\mu N)\}.
$$
The $\Ug$-action and the  ${\Bbb C}[{\frak h}^*]$-action on 
$\theta^+(\Res_\mu N)$  clearly induce similar actions on both 
$\widetilde{\theta}^+_\ell N$ and $\widetilde{\theta}^+_r N$.  Moreover, by 
construction, these actions fit together, making  $\widetilde{\theta}^+_\ell N$
and $\widetilde{\theta}^+_r N$  into $\widetilde{U}$-modules.  This way
we obtain two
functors  
$$\widetilde{\theta}^+_\ell,\, \widetilde{\theta}^+_r: \,\,
\Mod_\mu(\widetilde{U}) \to \Mod_\lambda(\widetilde{U})\,.$$
These functors are not
exact in general.  Notice however that if both $\lambda$  and $\mu$  are 
regular, then one has
$\widetilde{\theta}^+_r = \widetilde{\theta}^+_\ell = \theta^+$.

\proclaim{Proposition 1.5}  The functors  $\widetilde{\theta}^+_\ell$  and 
$\widetilde{\theta}^+_r$  are, respectively, the left and the right adjoint
of the functor  $\widetilde{\theta}^-$.  
\endproclaim

\demo{Proof}  Let  $\widetilde{M} \in \Mod_\lambda(\widetilde{U})$ and 
$\widetilde{N} \in \Mod_\mu(\widetilde{U})$.  Set  $M := \Res_\lambda \widetilde{M},$
$N := \Res_\mu \widetilde{N}$.  Then we have:
$$
\aligned
\Hom_{\Mod_\mu(\widetilde{U})} (\widetilde{\theta}^-\widetilde{M},\widetilde{N}) 
&= \{ f \in \Hom_{\Ug}(\theta^-M,N)\enspace \big\vert \enspace 
f(a*x) = a \cdot f(x),\quad
\forall a \in {\Bbb C}[{\frak h}^*]\} \\
&= \{ g \in \Hom_{\Ug}(M, \theta^+N)\enspace \big\vert\enspace
 g(a\cdot m) = a * g(m),\quad
\forall a \in {\Bbb C}[{\frak h}^*]\} .
\endaligned
$$
By  $\Ug$-linearity of the morphism  $g$  in the last line above,
and Lemma 1.2, we get 
$g(a\cdot m) = a \cdot g(m)$,\quad $\forall a \in {\Bbb C}[{\frak h}]^{W_\lambda}$.
Hence,  $a*g(m) = a\cdot g(m),\quad \forall a \in {\Bbb C}[{\frak h}^*]^{W_\lambda}$,
and  $\text{Image} (g) \subset \widetilde{\theta}^+_r N$.  Thus, $g \in 
\Hom_{\Mod_\lambda(\widetilde{U})} (\widetilde{M}, \widetilde{\theta}^+_r\widetilde{N})$.

The second adjunction is proved in a similar way.\quad $\square$
\enddemo
\vskip .2in


\subheading{2.\ Translation functors via ${\Cal D}$-modules}
\vskip .15in

For a complex variety $X$ let  ${\Cal O}_X$  and ${\Cal D}_X$  denote the 
sheaves of regular functions and regular differential operators on $X$ 
respectively.  

Let  $G \supset B \supset U$  be the  simply-connected semisimple Lie-group 
corresponding to ${\frak g}$,
a Borel subgroup of $G$,
 and the unipotent radical of $B$ respectively.  Let  $T = B/U$  be the abstract
maximal torus, see [CG, p.303].  
 Set  $\widetilde{\Cal B} = G/U$  and ${\Cal B} = G/B$, the 
flag manifold of $G$.  There is a natural $T$-action on $\widetilde{\Cal B}$
on the right that makes the canonical projection  $\pi: G/U \to G/B$  a 
principal $G$-equivariant $T$-bundle.

Let $\pi_{{}_\bullet} {\Cal D}_{\widetilde{\Cal B}}$  be the sheaf-theoretic 
direct image of ${\Cal D}_{\widetilde{\Cal B}}$ to ${\Cal B}$.  The right $T$-action 
on $\widetilde{\Cal B}$  induces a $T$-action on the sheaf 
$\pi_{{}_\bullet} {\Cal D}_{\widetilde{\Cal B}}$ by algebra automorphisms 
and we let  $\widetilde{\Cal D} \subset \pi_{{}_\bullet} {\Cal D}_{\widetilde{\Cal B}}$
denote the subsheaf of $T$-invariant sections
 of $\pi_{{}_\bullet} {\Cal D}_{\widetilde{\Cal B}}$.
Thus, $\widetilde{\Cal D}$  is a sheaf of algebras on ${\Cal B}$.

The infinitesimal left ${\frak g}$-action and the infinitesimal
right ${\frak h}$-action on 
$\widetilde{\Cal B}$  commute, giving rise to a homomorphism of the Lie 
algebra  ${\frak g} \times {\frak h}$  into global algebraic vector fields on 
$\widetilde{\Cal B}$.  This
 Lie algebra homomorphism can be extended to an associative
algebra homomorphism  $\Ug \otimes_{_{\Bbb C}} \Uh \to 
\Gamma(\widetilde{\Cal B},\,
{\Cal D}_{\widetilde{\Cal B}})^T$ (= {\it the algebra of right
$T$-invariant global
differential operators on}  $\widetilde{\Cal B}$).  It turns out that,
for any $z\in Z({\frak g})$,
the differential operator on ${\widetilde{\Cal B}}$
corresponding under the homomorphism above to the element
$z\otimes 1\in\Ug \otimes\Uh$ is equal
to the differential operator corresponding to the element
$1\otimes \chi(z)\in\Ug \otimes\Uh$, where $\chi: Z({\frak g})\iso
\Uh^W$ is the Harish-Chandra isomorphism. It follows that 
the homomorphism above factors through  $\Ug \otimes_{Z({\frak g})} \Uh 
\simeq \widetilde{U}$  (we identify  $\Uh$  with  ${\Bbb C}[{\frak
h}^*])$.
Furthermore,  the associated graded map of
"principal symbols" has been shown [BoBr]
to be a bijection, so that one obtains an algebra isomorphism:
$$
\widetilde{U}\ {\buildrel \sim\over\longrightarrow}\ 
\Gamma({\widetilde{\Cal B}},\,
{\Cal D}_{\widetilde{\Cal B}})^T = \Gamma({\Cal B}, \widetilde{\Cal D}).
\tag 2.1
$$

Observe that since the torus $T$ commutes with its own Lie algebra
action,
the image of $1\otimes \Uh\simeq {\Bbb C}[{\frak h}^*]$ is contained in the stalk of the
sheaf $\widetilde{\Cal D}$ at any point of ${\Cal B}$.
Notice further that the embedding ${\Bbb C}[{\frak h}^*] \hookrightarrow 
\widetilde{\Cal D}$ thus defined
 is central.  Hence, for any $\lambda \in {\frak h}^*$, we 
may define the category  $\Mod_\lambda(\widetilde{\Cal D})$  of coherent 
(sheaves on ${\Cal B}$ of) $\widetilde{\Cal D}$-modules $M$ such that 
$\J^n_\lambda \cdot M = 0$  for big enough  $n = n(M) \gg 0$.  Taking global 
sections of  $\widetilde{\Cal D}$-modules defines (via (2.1)) a functor 
$\Gamma_\lambda: \Mod_\lambda(\widetilde{\Cal D}) \to \Mod_\lambda(\widetilde{U})$.
Let 
$\Delta_\lambda: M \mapsto  \widetilde{\Cal D} \otimes_{\Gamma({\Cal B}, 
\widetilde{\Cal D})}M$ be the {\it localization functor}, which is
the left adjoint of $\Gamma_\lambda$.
We recall, see [BB] and also [BB3, Thm.3.3.1], the following important

\proclaim{Localization theorem 2.2}

\item{(i)} If $\lambda$ is $\rho$-dominant then the functor $\Gamma_\lambda$  is 
exact and  the canonical adjunction morphism
 $\Gamma_\lambda \cdot \Delta_\lambda \longrightarrow 
\id_{\Mod_\lambda(\widetilde{U})}$ is an isomorphism;

\item{(ii)} If $\lambda$  is regular then $\Gamma_\lambda$  gives an 
equivalence of the categories  $\Mod_\lambda(\widetilde{\Cal D})$ and 
$\Mod_\lambda(\widetilde{U})$, and the functor  $\Delta_\lambda$  is the 
inverse of  $\Gamma_\lambda$.  
\endproclaim

\demo{Remark 2.3}  To  $M \in \Mod_\lambda(\widetilde{\Cal D})$  assign 
a ${\Cal D}_{\widetilde{\Cal B}}$-module by the formula ${\Cal D}_{\widetilde{\Cal B}}
\otimes_{\widetilde{\Cal D}} M$.  This  ${\Cal D}_{\widetilde{\Cal B}}$-module
is clearly smooth along the fibres of the projection  $\pi: \widetilde{\Cal B} 
\to {\Cal B}$. That gives, for regular  $\lambda$, a fully faithful 
imbedding of the category  $\Mod_\lambda(\widetilde{U})$  into the category 
of $T$-monodromic ${\Cal D}_{\widetilde{\Cal B}}$-coherent modules. For
more details
see [BB3, $\S\S 2.5, 3.3$], and [Ka].\quad 
$\square$
\enddemo

\vskip 1pt

Recall that a subcategory ${\Cal A}$ of an abelian category ${\Cal C}$
is called a {\it
Serre subcategory} if ${\Cal A}$ is a full abelian subcategory stable
under taking extensions and subquotients in ${\Cal C}$. 
Given a Serre subcategory ${\Cal A}\subset {\Cal C}$, one can define
a quotient-category ${\Cal C}/{\Cal A}$. This is an abelian
category equipped  with an exact functor $\text{\it {quot}}: {\Cal C} \to
{\Cal C}/{\Cal A}$ such that ${\Cal A}=\,\text{Ker}\,(\text{\it {quot}}),$
where $\text{Ker}\,(\text{\it {quot}})$ stands for  the full subcategory of
${\Cal C}$ formed
by the objects $A$ such that $\text{\it {quot}}\,(A)=0$. 

Conversely, let $F: {\Cal C} \to {\Cal C}'$ be an exact functor between
abelian categories. Then, $\text{Ker}\,F$ is a Serre subcategory. Moreover,
there exists a unique functor $\overline{F} : {\Cal C}/\text{Ker}\,F \to
{\Cal C}'$,
such that the functor $F$ factors as the composition:
$$ {\Cal C} \ {\buildrel {\text{\it {quot}}}\over\longrightarrow}\
{\Cal C}/\text{Ker}\,F\ {\buildrel {\overline{F}}\over\longrightarrow}\
{\Cal C}'\,.$$

\proclaim{Lemma 2.4}  Let $F: {\Cal C} \to {\Cal C}'$ be an exact functor between
abelian categories which has a left (resp. right) adjoint functor
$F^\dag: {\Cal C}'\to{\Cal C}$. Then, $\overline{F} : {\Cal C}/\text{Ker}\,F \to
{\Cal C}'$ is an equivalence of categories if and only if the canonical
morphism $\id_{{\Cal C}'} \to F\cdot F^\dag$, (resp. 
$F^\dag\cdot F \to \id_{{\Cal C}}$) is an isomorphism.$\quad\square$
\endproclaim

From this lemma applied to the functor $F=\Gamma_\lambda$, and the
Localization
theorem 2.2, we deduce

\proclaim{Corollary 2.5}  Let $\lambda$ be $\rho$-dominant. Then, the
functor $\Gamma_\lambda$ induces an equivalence 
$$\overline{\Gamma}_\lambda :\,
\Mod_\lambda(\widetilde{\Cal D})/\text{Ker}\,\Gamma_\lambda
\iso \Mod_\lambda(\widetilde{U})\,.\quad\square$$
\endproclaim
\remark{Remark} See [Ka] for a more detailed description of the category
$\text{Ker}\,\Gamma_\lambda$ in the case of a non-regular $\lambda$.
$\quad\square$
\endremark

\vskip 3pt
Assume further
that  $\lambda \in {\frak h}^*$  is integral.  Then $\lambda$  gives rise 
to a homomorphism  $\dot{\lambda}: T \to {\Bbb C}^*$.  Let  ${\Cal O}(\lambda)$  
denote the sheaf on ${\Cal B}$  formed by all functions $f$  on $\widetilde{\Cal B}$  
such that:  
$$
f(x\cdot t) = \dot{\lambda}(t)f(x)\quad \forall x \in \widetilde{\Cal B},\ \ 
t \in T .   \tag 2.6
$$
Note that  ${\Cal O}(0) = {\Cal O}_{\Cal B}$, is
the structure sheaf of  ${\Cal B}$.  
More generally, ${\Cal O}(\lambda)$  is the sheaf of sections of a line bundle 
on ${\Cal B}$.  For a $\rho$-dominant $\lambda$,
 we have  $\Gamma({\Cal B}, {\Cal O}(\lambda)) = E_\lambda$,
an irreducible finite-dimensional ${\frak g}$-module with highest
weight  $\lambda$.

Let  ${\Cal D}_{\Cal B}(\lambda)$  denote the sheaf of 
differential operators acting on  ${\Cal O}(\lambda)$  (i.e., a sheaf of 
twisted differential operators on ${\Cal B}$).  Observe next that the sheaf 
${\Cal O}(\lambda)$  is stable under the natural action on functions of 
$T$-invariant differential operators on $\widetilde{\Cal B}$.  This gives an 
algebra morphism  $\widetilde{\Cal D} \longrightarrow {\Cal D}_{\Cal B}(\lambda)$.
One can show, see e.g., [BoBr],
that this morphism induces the following isomorphisms:
$$
{\Cal D}_{\Cal B}(\lambda) \simeq \widetilde{\Cal D}/
\widetilde{\Cal D}\!\cdot\!\J_\lambda
\quad \text{and}\quad \Gamma({\Cal B}, {\Cal D}_{\Cal B}(\lambda)) 
= \Ug/\Ug\!\cdot\!\I_{\la}\, .   \tag 2.7
$$

Given $\lambda, \mu \in {\frak h}^*$  such that  $\mu - \lambda$  is integral, 
define a geometric translation functor  $\Theta^\lambda_\mu : 
\Mod_\mu(\widetilde{\Cal D}) \to \Mod_\lambda(\widetilde{\Cal D})$  by 
the formula  $\Theta^\lambda_\mu: M \mapsto
{\Cal O}(\lambda - \mu) \otimes_{{\Cal O}_{\Cal B}}M$.
We have ${\Cal O}(\nu) \otimes_{{\Cal O}_{\Cal B}} {\Cal O} (-\nu) = 
{\Cal O}_{\Cal B}$,\ hence  $\Theta^\lambda_\mu \cdot \Theta^\mu_\lambda
= 
\id_{\Mod_\lambda(\widetilde{\Cal D})}$. 
Thus, the functor  $\Theta^\lambda_\mu$  is always an equivalence of categories.

Assume now that $\lambda$ and $\mu$  are integral $\rho$-dominant weights such 
that  $W_\lambda \subset W_\mu$.    Although the 
following result seems to be well known, we could not find its proof in the 
literature.

\proclaim{Proposition 2.8}   The following diagram commutes (up to
canonical equivalence of functors)
$$
\matrix \format \c &\quad \c &\quad \c \\
\Mod_\lambda(\widetilde{\Cal D}) & {\buildrel {\Gamma_\lambda}\over\longrightarrow} 
& \Mod_\lambda(\widetilde{U}) \\
\vspace{1\jot}
{\phantom{\Theta^\mu_\lambda}}\big\downarrow \Theta^\mu_\lambda && {\phantom{\Theta^\mu_\lambda}}\big\downarrow 
\widetilde{\theta}^\mu_\lambda \\
\vspace{1\jot}
\Mod_\mu(\widetilde{\Cal D}) & {\buildrel {\Gamma_\mu}\over\longrightarrow}
& \Mod_\mu(\widetilde{U}) 
\endmatrix
$$
\endproclaim

\demo{Proof} (essentially borrowed from [BB]):  Given a finite-dimensional 
${\frak g}$-module $E$, let  $E_{\Cal B} = E \otimes {\Cal O}_{\Cal B}$  be the trivial sheaf
of $E$-valued regular functions on ${\Cal B}$.  For such a function $f$ define 
a function  $\varphi_f: G \to E$  by the formula  $\varphi_f(g) = g^{-1}\cdot 
f(g)$.  The assignment  $f \longmapsto \varphi_f$ identifies  $E_{\Cal B}$ with 
the sheaf  $Ind^G_BE$  (on ${\Cal B}$) of germs of functions:
$$
\{ \varphi: G \to E \enspace\big\vert \enspace
\varphi(g\cdot b) = b^{-1} \cdot \varphi(g),\ \  
\forall b \in B\} .  \tag 2.9
$$
The infinitesimal action of ${\frak g}$ on $G$ by left translation makes  $Ind^G_BE$  a 
${\frak g}$-module.
\enddemo

By Lie's theorem, one
can find a $B$-stable filtration  $E_0 \subset E_1 \subset \ldots
\subset E_n = E$
 such that  $dim(E_i/E_{i-1}) = 1,$  for all $i$.  This gives a 
filtration on  $Ind^G_BE$  by the  ${\frak g}$-stable coherent subsheaves 
$Ind^G_BE_i$  defined by replacing  $E$ by $E_i$  in (2.9).  We have:
$$
Ind^G_BE_i/Ind^G_BE_{i-1} = Ind^G_B(E_i/E_{i-1}) = {\Cal O}(\nu_i)
$$
where  $\nu_i$  is the character of  $B$  corresponding to the 1-dimensional 
$B$-module  $E_i/E_{i-1}$.  

Assume now that  $E = E_{\mu -\lambda}$  is an irreducible ${\frak g}$-module 
with extreme weight $\mu -\lambda$  and  $M \in \Mod_\lambda(\widetilde{D})$.
We endow the sheaf  $(Ind^G_BE) \otimes_{{\Cal O}_{\Cal B}}M$  with the tensor 
product  ${\frak g}$-module structure.  The filtration  $E_0 \subset E_1 \subset
\ldots \subset E$  as above gives a ${\frak g}$-stable filtration on 
$(Ind^G_BE) \otimes_{{\Cal O}_{\Cal B}}M$  by the subsheaves\ 
$(Ind^G_BE_i) \otimes_{{\Cal O}_{\Cal B}}M$  with quotients of the form 
${\Cal O}(\nu) \otimes_{{\Cal O}_{\Cal B}}M$  where  $\nu$  is a weight of 
$E_{\mu -\lambda}$.  Clearly, for any  $a \in {\Bbb C}[{\frak h}^*]$,
for the action
of $a$ on ${\Cal O}(\nu) \otimes_{{\Cal O}_{\Cal B}}M$ we have the formula:
$$
a \cdot (f \otimes m) = f \otimes (T_\nu a) \cdot m,\qquad 
f \in {\Cal O}(\nu),\ \ m \in M\,,    \tag 2.10
$$
where $T_\nu$ denotes  the affine translation by $\nu$ on ${\Bbb
C}[{\frak h}^*]$ introduced
 before Proposition 1.4. In particular,  $Z({\frak g})$  acts on  
${\Cal O}(\nu) \otimes M$  via the (generalized) infinitesimal character\ 
$|\lambda + \nu| \in {\frak h}^*/W$.  

We now use the following result [BeGe, Lemma 1.5(iii)]: 

\noindent
{\it Let} $\nu, \nu' \in {\frak h}^*$  {\it be such that} $\lambda + \nu$ 
{\it is dominant}, $\lambda + \nu'$ {\it is}  $W$-{\it conjugate to}  
$\lambda + \nu$,  {\it and}  $\|\nu'\| \le \|\nu\|$  {\it in some euclidean}  
$W$-{\it invariant metric} $\| \cdot \|$  {\it on} ${\frak h}^*$.  {\it Then}  
$\nu' \in W_\lambda \cdot \nu$.  
\vskip 2pt

The result above implies that, for any weight  $\nu'$  of  $E_{\mu-\lambda}$  
other than  $\mu - \lambda$,  the point  $\lambda + \nu'$  is not $W$-conjugate 
to $\mu$.  Hence, the subquotient sheaf  ${\Cal O}(\mu - \lambda) 
\otimes_{{\Cal O}_{\Cal B}}M$  splits off from  $Ind^G_BE_{\mu -\lambda} 
\otimes_{{\Cal O}_{\Cal B}}M$  as a sheaf of  $Z({\frak g})$-modules, and moreover
we have:
$$
\text{pr}_{\ma}\bigl(Ind^G_B E_{\mu -\lambda} \otimes_{{\Cal O}_{\Cal B}}M\bigr) = 
{\Cal O}(\mu - \lambda) \otimes_{{\Cal O}_{\Cal B}}M,\qquad  M \in 
\Mod_\lambda(\widetilde{U})
\tag 2.11
$$
where  $\text{pr}_{\ma}$  stands for the projection to the $\ma$-isotypic
component of a $Z({\frak g})$-module.

Next, observe that, for any $E$ one has
$$
E \otimes \Gamma(M) = \Gamma(E_{\Cal B} \otimes_{{\Cal O}_{\Cal B}}M) = 
\Gamma(Ind^G_BE \otimes_{{\Cal O}_{\Cal B}}M).
$$
Thus, from (2.11) one obtains
$$
\aligned
\widetilde{\theta}^\mu_\lambda \cdot \Gamma_\lambda(M) &= 
\text{pr}_{\ma}\bigl(E_{\mu -\lambda} \otimes 
\Gamma(M)\bigr) \\
&= \text{pr}_{\ma}\Gamma\bigl(Ind^G_BE_{\mu -\lambda} \otimes_{{\Cal
O}_{\Cal B}}M\bigr)\\
&= \Gamma\bigl({\Cal O}(\mu - \lambda) \otimes_{{\Cal O}_{\Cal B}}M\bigr) = 
\Gamma_\mu \cdot \Theta^\mu_\lambda M.
\endaligned
$$
Finally, formula (2.10) shows that the ${\Bbb C}[{\frak h}^*]$-action on 
$\theta^\mu_\lambda \cdot \Gamma_\lambda(M)$  defined in section 1 corresponds to the 
natural ${\Bbb C}[{\frak h}^*]$-action on  $\Theta^\mu_\lambda M$.\quad $\square$

\vskip .2in


\subheading{3.\ Properties of translation functors}
\vskip .15in

Throughout this section we fix integral  $\rho$-dominant weights  $\lambda,\mu 
\in {\frak h}^*$  such that  $W_\lambda \subset W_\mu$,  and 
use the notation $\theta^+ := 
\theta^\lambda_\mu$  and $\theta^- := \theta^\mu_\lambda$, and similarly,
  $\Theta^+ := 
\Theta^\lambda_\mu$  and $\Theta^- := \Theta^\mu_\lambda$.

\demo{Proof of Proposition 1.4}  Let  $M \in \Mod_\lambda(\widetilde{U})$.
By Theorem 2.2(i), we can find a
$\widetilde{\Cal D}$-module  ${\Cal M}$  such that  $M =  \Gamma_\lambda({\Cal M})$.  
Part (ii) of the Proposition now follows from Proposition 2.8 and formula 
(2.10) for  $\nu = \mu - \lambda$.  
\enddemo

To prove part (i) take  $\widetilde{M} \in \Mod_\mu(\widetilde{U}),$
$\widetilde{N} \in \Mod_\lambda(\widetilde{U})$  and set  $M := 
\Res_\mu \widetilde{M},$ $N := \Res_\lambda \widetilde{N}$.  An element 
$a \in {\Bbb C}[{\frak h}]^{W_\mu}$  gives rise to an endomorphism  
$\theta^+(a) \in$
$ \text{End}\,(\theta^+\!M)$, hence, to an endomorphism  $\theta^+(a)_{{}_M}$ 
of  $\Hom_{\Ug}(\theta^+M, N)$, via composition with $\theta^+(a)$.  We have a 
commutative diagram:
$$
\matrix \format \c & \c & \c & \c & \c & \c 
& \c & \c & \c\\
\Hom(\theta^+M, N) & = & \Hom(M, \theta^- N) &=& \Hom(M, \theta^-N) &=& 
\Hom(M, \theta^-N)  &=& \Hom(\theta^+ M, N) \\
\vspace{1\jot}
{\phantom{a_M}}\big\downarrow \theta^+(a)_{{}_M} && 
{\phantom{a_{{}_M}}}\big\downarrow a_{{}_M} && {\phantom{a_M}}\big\downarrow 
a_{{}_{\theta^{^-}\!N}} && {\phantom{aM}}\big\downarrow 
\theta^-(T_{\mu -\lambda}a) && {\phantom{aM}}\big\downarrow 
(T_{\mu -\lambda}a)_{{}_{\theta^{^+}\!M}}  \\
\vspace{1\jot}
\Hom(\theta^+M, N) & = & \Hom(M, \theta^-N) &=& \Hom(M, \theta^-N) &=& 
\Hom(M, \theta^-N)  &=& \Hom(\theta^+ M, N) .
\endmatrix
$$
Thus, we deduce that
$\theta^+(a) = T_{\mu - \lambda}a$.\quad  $\square$
\medpagebreak

Here are the most important properties
of the translation functors.

\proclaim{Proposition  3.1}  (i)
The functor 
$\widetilde{\theta}^- :\,\Mod_\lambda(\widetilde{U})\to
\Mod_\mu(\widetilde{U})$ induces an equivalence:$\qquad$
$\Mod_\lambda(\widetilde{U})/\text{Ker}\,\widetilde{\theta}^-
\iso\Mod_\mu(\widetilde{U})$;

(ii) The adjunction morphisms induce  isomorphisms of functors:
$$\widetilde{\theta}^- \cdot 
\widetilde{\theta}^+_\ell = \id_{\Mod_\mu(\widetilde{U})}
=\widetilde{\theta}^- \cdot 
\widetilde{\theta}^+_r\,;$$

(iii) We have a natural isomorphism of functors:
$\quad\widetilde{\theta}^+_\ell = \Gamma_\lambda \cdot 
\Theta^+ \cdot \Delta_\mu$.
\endproclaim

\demo{Proof} We begin with a general remark. Let ${\Cal C}$ be 
an abelian category, and
${\Cal C}_2 \subset {\Cal C}_1 \subset {\Cal C}$ two Serre
subcategories. Then the various  quotient categories are related
by the following canonical transitivity isomorphism
$$ {\Cal C}/{\Cal C}_1 \;\simeq\; 
({\Cal C}/{\Cal C}_2)/({\Cal C}_1/{\Cal C}_2)\,.\tag 3.1.1$$

We set ${\Cal C}:=\Mod_\lambda(\widetilde{\Cal D})$.
Since the geometric translation functor
$\Theta^-:\,\Mod_\lambda(\widetilde{\Cal D})$
$\iso
\Mod_\mu(\widetilde{\Cal D})$ is an equivalence of categories,
we may (and will) identify the category $\Mod_\mu(\widetilde{\Cal D})$
with ${\Cal C}$.
Further, by Corollary 2.5, we have equivalences
$$
\Mod_\lambda(\widetilde{\Cal D})/\text{Ker}\,\Gamma_\lambda
\iso \Mod_\lambda(\widetilde{U})\quad,\quad
\Mod_\mu(\widetilde{\Cal D})/\text{Ker}\,\Gamma_\mu
\iso \Mod_\mu(\widetilde{U})\,.$$
Thus,  setting  ${\Cal C}_1:=\text{Ker}\,\Gamma_\mu$ and
${\Cal C}_2:=\text{Ker}\,\Gamma_\lambda$ we can rewrite the above as
follows
$${\Cal C}/{\Cal C}_1 \simeq\Mod_\mu(\widetilde{U})\quad,\quad
{\Cal C}/{\Cal C}_2\simeq\Mod_\lambda(\widetilde{U})\,.\tag 3.1.2$$
Moreover, Proposition 2.8 insures that 
${\Cal C}_2 \subset {\Cal C}_1$ and that the natural functor
$ {\Cal C}/{\Cal C}_1 \to {\Cal C}/{\Cal C}_2$ gets identified under
the equivalences (3.1.2) with the functor $\widetilde{\theta}^-:
\Mod_\lambda(\widetilde{U}) \to \Mod_\mu(\widetilde{U})$.
Part (i) of Proposition 3.1 now follows from the
transitivity equivalence (3.1.1) applied to our categories.
With part (i) being established, Lemma 2.4  yields
part (ii) of Proposition 3.1.

To prove (iii) we take
adjoints in the diagram of  Proposition 2.8. This
yields the following commutative diagram:
$$
\matrix \format \c &\quad \c &\quad \c \\
\Mod_\lambda(\widetilde{\Cal D}) & {\buildrel {\Delta_\lambda}\over
\longleftarrow} & \Mod_\lambda(\widetilde{U}) \\
\vspace{1\jot}
\Theta^+\big\uparrow {\phantom{\Theta^+}}&& {\phantom{\Theta^+}}\big\uparrow 
\widetilde{\theta}^+_\ell \\
\vspace{1\jot}
\Mod_\mu(\widetilde{\Cal D}) & {\buildrel {\Delta_\mu}\over
\longleftarrow} & \Mod_\mu(\widetilde{U}) . 
\endmatrix  \tag 3.2
$$
The diagram implies readily:
 $\;\widetilde{\theta}^+_\ell = \Gamma_\lambda \cdot 
\Theta^+ \cdot \Delta_\mu\,,$ where we have used that every object of
the category $\Mod_\mu(\widetilde{\Cal D})$ is isomorphic to an object
of the form $\Delta_\mu(M)$, by Theorem 2.2(i). Part (iii) of
Proposition 2.8 is proved.

One may reverse the logic slightly to derive parts (i) and (ii)
of Proposition 3.1 in an alternative way, as follows.
First of all one establishes part (iii) of the
proposition the same way as above. With the formula of part (iii)
at hand we calculate:
$$\widetilde{\theta}^- \cdot
\widetilde{\theta}^+_\ell \;=\; \widetilde{\theta}^- \cdot \Gamma_\lambda \cdot
\Theta^+ \cdot \Delta_\mu \;=\; \Gamma_\mu \cdot \Theta^- \cdot \Theta^+ \cdot 
\Delta_\mu \;=\; \Gamma_\mu \cdot \Delta_\mu \;=\; \id_{\Mod_\mu(\widetilde{U})}\,,$$ 
where the first equality holds by Proposition 3.1(iii),
the second is Proposition 2.8, and the last one is
due to Theorem 2.2(i). This gives
part (ii), and part (i) now follows from the "if" part of Lemma 2.4.
\quad $\square$
\enddemo

\proclaim{Corollary 3.3}  For  $M \in \Mod_\lambda(\widetilde{U})$,  the 
following conditions are equivalent:

\item{(i)}  For any  $N \in \Mod_\lambda(\widetilde{U})$, the following natural 
morphism is an isomorphism
$$
\Hom_{\Mod_\lambda(\widetilde{U})}(M,N) \iso
\Hom_{\Mod_\mu(\widetilde{U})}(\widetilde{\theta}^-M, \widetilde{\theta}^-N)\;;
$$

\item{(ii)}  There exists an $M' \in \Mod_\mu(\widetilde{U})$  such that
$M = \widetilde{\theta}^+_\ell(M')$.  
\endproclaim

\demo{Proof}  Observe that 
the canonical
adjunction morphism  $\widetilde{\theta}^+_\ell \cdot 
\widetilde{\theta}^-M \to M$  is an isomorphism if and only if the induced 
morphism
$$
\Hom(M,N) \to \Hom(\widetilde{\theta}^+_\ell \cdot \widetilde{\theta}^-M, N) 
= \Hom(\widetilde{\theta}^-M, \widetilde{\theta}^-N)
$$
is an isomorphism, for every $N$.  In the latter
 case we may put  $M' = \widetilde{\theta}^-M$.  
Conversely, if  $M = \widetilde{\theta}^+_\ell(M')$  then we have 
$\widetilde{\theta}^+_\ell \cdot \widetilde{\theta}^-M = \widetilde{\theta}^+_\ell
\cdot (\widetilde{\theta}^- \cdot \widetilde{\theta}^+_\ell M') = 
\widetilde{\theta}^+_\ell M' = M$  (by Proposition 3.1(ii)),  and the result 
follows.\quad $\square$
\enddemo

\medpagebreak

Recall that to any  $M \in \Mod_{\ma}(\Ug)$  we can associate (cf. Lemma 1.3)
the $\widetilde{U}$-module 
$\widetilde{M} := \widetilde{U} \otimes_{\widetilde{U}^{W_\mu}}M$.
Note that the
module $\widetilde{M}$ has
a natural $W_\lambda$-action induced from the  $W_\lambda$-action on $\widetilde{U}$, 
which commutes with the  $\Ug$-action.  By functoriality, this gives 
a $W_\lambda$-action on the $\Ug$-module  $\Res_\lambda \cdot 
\widetilde{\theta}^+_\ell(\widetilde{M})$. 
We have:

\proclaim{Proposition 3.4}  There is a functorial isomorphism:
$$
\theta^+M = \Bigl(\Res_\lambda \cdot\, \widetilde{\theta}^+_\ell(\widetilde{U} 
\otimes_{\widetilde{U}^{W_\mu}}M)\Bigr)^{W_\lambda},\quad  M \in 
\Mod_{\ma}({U}).
$$
\endproclaim

\demo{Proof}  It is clear that there is
 a functor isomorphism $\theta^- \cdot \Res_\lambda = 
\Res_\mu \cdot\, \widetilde{\theta}^-:  \Mod_\lambda(\widetilde{U})\iso
\Mod_\ma(\Ug)$.  Taking adjoints on each side,
and using Lemma 1.3  we obtain:
$$\widetilde{\theta^+M} =\widetilde{U} 
\otimes_{\widetilde{U}^{W_\lambda}} (\theta^+M) =
\widetilde{\theta}^+_\ell(\widetilde{U} 
\otimes_{\widetilde{U}^{W_\mu}}M)=
\widetilde{\theta}^+_\ell(\widetilde{M})\,.$$
The isomorphisms above are  compatible with the $W_\lambda$-actions.  
We take $W_\lambda$-invariants on each side of the isomorphism. 
Note that
$\widetilde{U}$  is a free $\widetilde{U}^{W_\lambda}$-module
isomorphic
to $\C[W_\lambda]\otimes_{_\C}\widetilde{U}^{W_\lambda}$
as a $W_\lambda$-module. Hence
taking $W_\lambda$-invariants of the leftmost term in  the above isomorphisms
yields
$\bigl(\widetilde{U} 
\otimes_{\widetilde{U}^{W_\lambda}} (\theta^+M)\bigr)^{W_\lambda}=
\theta^+M$. Comparing with the $W_\lambda$-invariants of
the rightmost term
completes the proof.\quad $\square$
\enddemo

\medpagebreak

\proclaim{Proposition 3.5}  There is a functorial isomorphism
$$
\theta^- \cdot\, \theta^+M \simeq \Res_\mu(\widetilde{U}^{W_\lambda} 
\otimes_{\widetilde{U}^{W_\mu}}M).
$$
\endproclaim

\demo{Proof}  As in the proof of Proposition 3.4 we have 
$\widetilde{\theta^+M} = 
\widetilde{\theta}^+_\ell (\widetilde{M})$.
The functor  $\theta^-$  clearly commutes with the functor  $\Res \cdot 
(\widetilde{U} \otimes_{\widetilde{U}^{W_\mu}} \bullet)$ and with the 
$W$-actions.  Hence, we obtain 
$$
\aligned
\Res_\mu \Bigl(\widetilde{U} \otimes_{\widetilde{U}^{W_\lambda}}(\theta^-\!\cdot
\theta^+M)\Bigr) &= \theta^- \cdot \Res_\lambda \cdot\, 
\widetilde{\theta}^+_\ell(\widetilde{M}) \\
&= \Res_\mu \cdot \,\widetilde{\theta}^- \cdot 
\widetilde{\theta}^+_\ell(\widetilde{M}) = 
\Res_\mu(\widetilde{M}).
\endaligned
$$
Taking  $W_\lambda$-invariants on each side completes the proof.
\quad $\square$\enddemo

\vskip .2in

\subheading{4.\ Applications to the category ${\Cal O}$}
\vskip .1in

Fix a Cartan and Borel subalgebras ${\frak h}\subset
{\frak b}\subset {\frak g}$. Thus we may identify this Cartan subalgebra
with the abstract Cartan subalgebra ${\frak b}/[{\frak b},{\frak b}]$
via the composition ${\frak h}\hookrightarrow
{\frak b}\twoheadrightarrow
{\frak b}/[{\frak b},{\frak b}]\,.$
Given $\lambda \in {\frak h}^*$, we define the category ${\Cal
O}_\lambda$
as the full subcategory of $\Mod_\la(\Ug)$ formed by the $\Ug$-modules
$M$
such that

\noindent
$\bullet\quad$ $\text{U}\!{\frak b}$-action on $M$ is locally finite, and

\noindent
$\bullet\quad$ $\text{U}\!{\frak h}$-action on $M$ is diagonalizable.

Given a weight
$\mu\in {\frak h}^*$, view it as 
a 1-dimensional ${\frak b}$-module, ${\Bbb C}_\mu$,
via the projection ${\frak b}\twoheadrightarrow
{\frak b}/[{\frak b},{\frak b}]=
{\frak h}\,{\buildrel\mu\over
\longrightarrow}\,{\Bbb C}\,.$
Let $M_\mu :=\Ug \otimes_{\text{U}\!{\frak b}}{\Bbb C}_\mu$ 
be the Verma module with highest weight
$\mu$. We write $L_\mu$ for its simple quotient.
Then, for any $w\in W$, we have $M_{w\cdot\lambda}\,,\,L_{w\cdot\lambda}\in
{\Cal O}_\lambda$. The category 
${\Cal O}_\lambda$  is known [BGG] to have enough projectives, and we
let
 $P_{w\cdot\lambda}$ denote the
indecomposable projective cover of $L_{w\cdot\lambda}$ in 
${\Cal O}_\lambda$. 
We refer to [BGG] for more properties of
the category ${\Cal
O}_\lambda$.

\medskip

Clearly  $-\rho$  is the fixed point of the dot-action.  The Verma module 
$M_{-\rho}$  is simple and is the unique simple object of the category 
${\Cal O}_{-\rho}$.  Moreover, $M_{-\rho}$  is also a projective in 
${\Cal O}_{-\rho}$, for $-\rho$  is a $\rho$-dominant weight.  Thus,
${\Cal O}_{-\rho}$  is a semisimple category and any object of 
${\Cal O}_{-\rho}$  is isomorphic to a direct sum of finitely many
copies of 
$M_{-\rho}$.  

Now let $\lambda$ and $\mu$  denote integral $\rho$-dominant weights such 
that  $W_\lambda \subset W_\mu$.  The adjoint functors $(\theta^\mu_\lambda, 
\theta^\lambda_\mu): {\Cal O}_\lambda \rightleftarrows {\Cal O}_\mu$  are exact,
hence, take projectives into projectives.  Further, it is known  and easy
to prove, see e.g. [Ja],  that, for any 
$w \in W$, one has $\theta^\mu_\lambda M_{w \cdot \lambda} = M_{w \cdot
\mu}$.
Therefore, by 
the exactness, one gets:
$$
\theta^\mu_\lambda L_{w \cdot\lambda} = \cases L_{w \cdot \mu} &\quad \text{if}\quad 
w \in W_\mu \cdot w_0 \\
0 &\quad \text{if}\quad 
w \ne W_\mu \cdot w_0 .
\endcases
$$
It follows, by adjunction, that the exact functor  
$\Hom_{{\Cal O}_\lambda}(\theta^\lambda_\mu P_\mu,\bullet \ )$ on ${\Cal O}_\lambda$ 
kills all simple modules but  $L_\lambda$.  Thus, we obtain

\proclaim{Proposition 4.1}  $\theta^\lambda_\mu P_\mu = P_\lambda$.  
\endproclaim

From now on let $\lambda$ be an integral $\rho$-dominant weight.  Proposition 
4.1 yields (for $\mu = - \rho)$:
$$
P_\lambda = \theta^\lambda_{-\rho} M_{- \rho}.  \tag 4.2
$$

Let  $I^W$  be the ideal in  ${\Bbb C}[{\frak h}^*]$  generated by all 
$W$-invariant polynomials without constant term.  Set  $C = {\Bbb C}[{\frak h}^*]/I^W$,
the coinvariant algebra.  The $W$-action on ${\Bbb C}[{\frak h}^*]$  induces a 
$W$-action on $C$  making it a regular representation of $W$.  For each 
$\lambda \in {\frak h}^*$  we have the subalgebra  $C^{W_\lambda} \subset C$.

The results of the previous section enable us to give a short proof of the 
following important theorem that was implicitly conjectured in [BG] and
was first 
proved in [S1, Thm. 3] (see also [Be]).  

\proclaim{Endomorphism-theorem 4.3}  There is a canonical algebra isomorphism
$$
End_{\Ug} P_\lambda \simeq C^{W_\lambda}.
$$
\endproclaim

\demo{Proof}  Set  $\theta^+ = \theta^\lambda_{-\rho}$  and  $\theta^- = 
\theta_\lambda^{-\rho}$.  From (4.2) and Proposition 3.5 we get:
$$
\aligned
\Hom_{\Ug}(P_\lambda, P_\lambda) &= \Hom_{\Ug}(\theta^+M_{-\rho}, 
\,\theta^+M_{-\rho}) =\\
&= \Hom_{\Ug}(M_{-\rho},\, \theta^- \cdot \theta^+ M_{-\rho})\\
&= \Hom_{\Ug}\Bigl(M_{-\rho},\ \Res_{-\rho}(\widetilde{U}^{W_\lambda} 
\otimes_{\Ug} M_{-\rho})\Bigr) \\
&=
\widetilde{U}^{W_\lambda}/
\widetilde{U}^{W_\lambda}\!\cdot\!\J_{-\rho} 
\simeq C^{W_\lambda}.
\endaligned
$$

To complete the proof we must
 show that the chain of isomorphisms above gives rise to an {\it algebra} map
$\End_{\Ug} P_\lambda$
$ \to C^{W_\lambda}$. To  this end, observe that
the functor $\theta^-$ induces a ring homomomorphism:
$$\End_{\Ug} P_\lambda \to
\End_{\widetilde{U}^{W_\lambda}}(\theta^-P_\lambda)
=\End_{\widetilde{U}^{W_\lambda}}(\theta^-\theta^+M_{-\rho})=
\End_{\widetilde{U}^{W_\lambda}}(\widetilde{U}^{W_\lambda} 
\otimes_{\Ug} M_{-\rho})\,.\tag 4.3.1$$

Observe further that
the action of the
central
subalgebra $S{\frak h}^{W_\lambda}\subset \widetilde{U}^{W_\lambda}$ on 
$\widetilde{U}^{W_\lambda} 
\otimes_{\Ug} M_{-\rho}$  gives an algebra isomorphism
$\tau: C^{W_\lambda} \iso \End_{\widetilde{U}^{W_\lambda}}(\widetilde{U}^{W_\lambda} 
\otimes_{\Ug} M_{-\rho})$,
since $M_{-\rho}$ is a simple $\Ug$-module.
The result now follows by composing  (4.3.1) with $\tau^{-1}.\quad\square$
\enddemo
\vskip 3pt

Let $B\subset G$ be the Borel subgroup corresponding to the Borel
subalgebra ${\frak b}$.
Write ${\Cal B}_w = B \cdot w \cdot B/B \subset {\Cal B}$  for the Bruhat cell in 
the Flag manifold corresponding to  $w \in W$. Let  $j_w:  {\Cal B}_w \hookrightarrow
{\Cal B}$  denote the embedding, and ${\Cal M}^*_w := (j_w)_* 
{\Cal O}_{{\Cal B}_w}$, the standard ${\Cal D}_{{\Cal B}}$-module ``of 
holomorphic distributions", supported on $\overline{{\Cal B}}_w$.  Let  
$\DD M_{w \cdot \lambda}$  denote the dual (in the category 
${\Cal O}$) of the Verma module $M_{w\cdot \lambda}$.  

The following result is known, see e.g. [BB3].

\proclaim{Proposition 4.4}  Let $\lambda$  be an integral $\rho$-dominant 
weight.  Then, for any  $w \in W^\lambda$  and $y \in W$,  there are natural 
isomorphisms of  ${\Cal D}(\lambda)$- and $\Ug$-modules respectively:
$$
\Delta_\lambda (\DD M_{w \cdot \lambda}) \simeq {\Cal O}(\lambda) 
\otimes_{{\Cal O}_{\Cal B}} {\Cal M}^*_{w w_0}\quad \text{$\text{and}$}\quad
\Gamma({\Cal B},\,{\Cal O}(\lambda) \otimes {\Cal M}^*_{y w_0})\, \simeq\, 
\DD M_{y \cdot \lambda}\,.\enspace\square
$$
\endproclaim

\bigpagebreak

Following [S1], define an exact functor  ${\Bbb V}: {\Cal O}_\lambda 
\to C^{W_\lambda}$-$\Mod$  by the formula  ${\Bbb V}: M \longmapsto 
\Hom_{\Ug}(P_\lambda, M)$, where the  $Hom$-space is viewed as an 
$End_{\Ug}P_\lambda$-module, hence a $C^{W_\lambda}$-module, via 
composition.  

One may reinterpret the functor ${\Bbb V}$  as follows.  We have by 
(4.2), ${\Bbb V}(M) =$  

\noindent
$\Hom_{\Ug}(M_{-\rho}, \theta_\lambda^{-\rho}M)$.  
The  $C^{W_\lambda}$-module structure on  ${\Bbb V}(M)$  arises from a natural 
${\Bbb C}[{\frak h}^*]^{W_\lambda}$-action on $\theta_\lambda^{-\rho}M$.  The 
latter action was actually defined in n.1.  Namely, view $M$ as a 
$\widetilde{U}^{W_\lambda}$-module, via Lemma 1.2.  This gives a 
${\Bbb C}[{\frak h}]^{W_\lambda}$-action on $M$ which induces a 
${\Bbb C}[{\frak h}^*]^{W_\lambda}$-action on  $\theta_\lambda^{-\rho}M$, 
by functoriality.  Furthermore, the subalgebra  $Z({\frak g}) \simeq 
{\Bbb C}[{\frak h}]^W \subset {\Bbb C}[{\frak h}^*]^{W_\lambda}$  acts 
trivially on  $\Hom_{\Ug}(M_{-\rho}, \theta_\lambda^{-\rho}M)$, 
for it obviously acts trivially on  $M_{-\rho}$.  Thus, the action on 
${\Bbb V}(M)$  factors through  $C^{W_\lambda}$.

Using the results of the previous sections we can simplify the proof of the 
key theorem [S1, Struktursatz 2] saying that the functor  ${\Bbb V}$  is 
faithful on projectives and injectives; we have 

\proclaim{Theorem 4.5}  For any injective module  $I \in {\Cal O}_\lambda$ 
and any  $M \in {\Cal O}_\lambda$, the natural morphism
$$
\Hom_{\Ug}(I, M) \longrightarrow \Hom_{C^{W_\lambda}}({\Bbb V}(I), 
{\Bbb V}(M))
$$
is an isomorphism.
\endproclaim

\demo{Remark 4.6}  The functor  $\theta_\lambda^{-\rho}$, hence ${\Bbb V}$,
clearly commutes with the standard duality  $N \longmapsto \DD N$ 
on the category ${\Cal O}$.  Furthermore, $I$ is an injective in 
${\Cal O}_\lambda$ if and only if  $\DD I$  is a projective in 
${\Cal O}_\lambda$.  Thus, dualizing Theorem 4.5 we get that, for any 
projective  $P \in {\Cal O}_\lambda$  and any  $M \in {\Cal O}_\lambda$, the 
natural morphism
$$
\Hom_{\Ug}(M,P) \longrightarrow \Hom_{C^W_\lambda}({\Bbb V}(M), 
{\Bbb V}(P))
$$
is an isomorphism.$\quad\square$
\enddemo

\medpagebreak

Let  $M \mapsto  \widetilde{M} = \widetilde{U} \otimes_{\widetilde{U}^{W_\lambda}}M$
be the left adjoint of the functor  $\Res_\lambda: \Mod_\lambda(\widetilde{U}) 
\to \Mod_{\la}(\Ug)$, cf. Lemma 1.3.  The $W_\lambda$-action 
on $\widetilde{U}$  induces a $W_\lambda$-action on $\widetilde{M}$ and, for 
any  $M,N \in \Mod_\la(\Ug)$, there is a canonical isomorphism
$$
\Hom_{\Ug}(M,N)\ {\buildrel \sim\over\longrightarrow}\ 
\Hom_{\widetilde U}(\widetilde{M}, \widetilde{N})^{W_\lambda}.  \tag 4.7
$$

We will deduce Theorem 4.5 from the following result.

\proclaim{Lemma 4.8}  For any injective  $I \in {\Cal O}_\lambda$, the 
adjunction morphism  $\widetilde{\theta}_\ell^+ \cdot 
\widetilde{\theta}^-(\widetilde{I}) \to \widetilde{I}$\ 
(where  $\theta^+ = \theta^\lambda_{-\rho}$)  is an isomorphism.
\endproclaim

\demo{Proof of the Lemma}  Let  $\DD M_\lambda$  be the dominant dual 
Verma module, viewed as a $\widetilde{U}$-module via the projection 
$\widetilde{U} \to \widetilde{U}/\widetilde{U}\!\cdot\! \J_\lambda \simeq
\Ug/\Ug \!\cdot\! \J_\la$, and let 
$M_{-\rho}$  be viewed as a $\widetilde{U}$-module in a similar way.
By Proposition 3.1(iii) we have, $\widetilde{\theta}^+_\ell M_{-\rho} = 
\widetilde{\theta}^+_\ell(\DD M_{-\rho}) 
= \Gamma_\lambda \cdot \Theta^\lambda_{-\rho} \cdot
\Delta_{-\rho}(\DD M_{-\rho})$.  Hence,    Proposition 4.4 yields
$$
\widetilde{\theta}^+_\ell(M_{-\rho}) = \DD M_\lambda.
$$
It follows from Corollary 3.3 that the adjunction morphism 
$\widetilde{\theta}^+_\ell \cdot \widetilde{\theta}^-(\DD M_\lambda) \to 
\DD M_\lambda$  is an isomorphism.

Now, $\DD M_\lambda$  is an injective in ${\Cal O}_\lambda$.
Furthermore,  it was 
shown in [BeGe, Thm. 3.3(b)] that any injective in ${\Cal O}_\lambda$  can be written in 
the form  $I = \Phi(\DD M_\lambda)$  for an appropriate projective 
functor $\Phi$ (see $\S 6$ below).
  Hence, $\widetilde{I} = \widetilde{\Phi}(\DD M_\lambda)$ 
and the lemma follows from Proposition 6.6 (of $\S 6$),
which is independent of the intervening material. \quad $\square$
\enddemo

\demo{Proof of Theorem 4.5}  Let  $I$  be an injective in ${\Cal O}_\lambda$.
Then, for any $M \in {\Cal O}_\lambda$, from Lemma 4.8 and Corollary 3.3 we 
obtain an isomorphism
$$
\Hom_{\widetilde U}(\widetilde{I}, \widetilde{M})\ {\buildrel \sim\over\longrightarrow}\ 
\Hom_C(\widetilde{\theta}^-(\widetilde{I}), \,\widetilde{\theta}^-(\widetilde{M})).
$$
Taking $W_\lambda$-invariants on each side and using (4.7) completes the 
proof.\quad $\square$
\enddemo
\vskip .2in

\subheading{5.\ Convolution of Harish-Chandra modules}
\vskip .1in

Throughout this section $\lambda$ and $\mu$  stand for integral $\rho$-dominant 
weights.

\medpagebreak

Let  $\Ug_{_{\la}}$ and  $\widetilde{U}_\lambda$  denote the 
completions of the algebras  $\Ug$  and $\widetilde{U}$  with 
respect to the  $\I_{\la}$-adic and  $\J_\lambda$-adic 
topology, respectively
(recall that $\I_{\la}\subset Z({\frak g})$
and $\J_\lambda\subset{\Bbb C}[{\frak h}^*]\,$ are maximal ideals).
Since $\Ug$ is a Noetherian algebra, its $\I_{\la}$-adic completion,
$\Ug_{_{\la}}$ is also Noetherian. 
Hence, the Artin-Rees lemma implies that any finitely-generated
$\Ug_{_{\la}}$-module is complete with respect to the $\I_{\la}$-adic
topology. Therefore, finitely-generated
$\Ug_{_{\la}}$-modules form an abelian category. Similar considerations
apply to the algebra $\widetilde{U}_\lambda$.

 By a {\it complete}  
$\Ug$-module, resp. $\widetilde{U}$-module,
we will mean an $\Ug$-module $M$ which is complete
in the $\I_{\la}$-adic, resp. $\J_\lambda$-adic,
topology, i.e., such that \hfill$M\;=$\newline\noindent
$\lim\text{proj}\;~(M/\I_{_\la}^n\!\cdot\!M)$. A complete $\Ug$-module
  is not necessarily 
annihilated by some power of the ideal  $\I_{\la}$,
but it is isomorphic to 
a limit of the projective system of  $\Ug$-bimodules:  
$M/\I^1_{\la}\!\cdot\!M \leftarrow$
$M/\I^2_{\la}\!\cdot\!M\leftarrow \ldots$.
The algebra $\Ug_{_{\la}}$ is an example
of a complete $\Ug$-module. 
Furthermore, any complete 
 $\Ug$-module $M$ has a natural structure of a $\Ug_{_{\la}}$-module,
and we will often make no distinction between these $\Ug$-module
and $\Ug_{_{\la}}$-module structures on  $M$.

Below, finitely generated
left $\Ug_{_{\la}}\otimes_{\Bbb C}\Ug_{_{\la}}^{opp}$-modules will
be referred to as finitely generated
$\Ug_{_{\la}}$-bimodules, and similar
terminology will be used for $\widetilde{U}_\lambda
\otimes_{\Bbb C}\widetilde{U}_\lambda^{opp}$-modules.
Let  $\HC_{\la}$, resp. $\widetilde{\HC}_{\lambda}$,
be the category of finitely-generated {\it complete}  
$\Ug_{_{\la}}$-bimodules
(resp. $\widetilde{U}_\lambda$-bimodules) 
$M$ such that the adjoint ${\frak g}$-action 
$\ad\ x : m \longmapsto x \cdot m - m \cdot x$  on $M$  is locally-finite.  Any 
object of  $\HC_{\la}$  is 
finitely-generated both as a left 
and as a right  $\Ug_{_{\la}}$-module.
  Therefore, the Artin-Rees
lemma insures,
as explained in the first paragraph of this section, 
that $\HC_{\la}$ is an abelian category.
 Similar remarks apply to the category 
$\widetilde{\HC}_{\lambda}$.

\medpagebreak

To localize complete
$\Ug_{_{\la}}$-bimodules we recall the imbedding
${\Bbb C}[{\frak h}^*] \subset \widetilde{\Cal D}$, and
form the sheaf 
$\widetilde{\Cal D}_\lambda$, the $\J_\lambda$-adic completion of the algebra 
$\widetilde{\Cal D}$. There is an algebra isomorphism
$\Gamma(\Cal{B},\,\widetilde{\Cal D}_\lambda) 
= \widetilde{U}_\lambda$.  For any sheaf ${\Cal M}$ 
 of coherent  $\widetilde{\Cal D} 
\boxtimes \widetilde{\Cal D}^{opp}$-modules on ${\Cal B} \times {\Cal B}$, the 
space $\Gamma({\Cal B} \times {\Cal B},\,
{\Cal M})$  has a natural $\widetilde{U}\otimes
\widetilde{U}^{opp}$-module structure.  
However, the assignment  ${\Cal M} \mapsto  \Gamma({\Cal B} \times {\Cal
B},\,
{\Cal M})$  is not an exact 
functor, for $\Gamma$  is not exact on {\it right} $\widetilde{\Cal D}$-modules.
Thus, one has to localize $\widetilde{U}$-bimodules in a different way,
as we now explain.

We begin with a well-known observation that if ${\Cal M}$ is a left
${\Cal D}_X$-module on an algebraic variety $X$, and $\Omega_X$ is the
line bundle of top-degree regular forms on $X$, then the sheaf
$\Omega_X\otimes_{{\Cal O}_X}{\Cal M}$ has a natural right
${\Cal D}_X$-module structure. Put another way, there is a canonical
algebra isomorphism 
$${\Cal D}_X^{opp} \iso\Omega_X
\otimes_{{\Cal O}_X}{\Cal D}_X\otimes_{{\Cal O}_X}\Omega_X^{-1}\,.\tag 5.1$$
Now, let $X={\Cal B}$ be the flag manifold. Then $\Omega_{{\Cal B}}=
{\Cal O}(-2\rho)$; further, it is  clear  that, for any integral
$\lambda$, there is a natural isomorphism
$$\widetilde{\Cal D}_\lambda  \simeq {\Cal O}(\lambda)
\otimes_{{\Cal O}_{{\Cal B}}}
{\Cal D}_{{\Cal B}}\otimes_{{\Cal O}_{{\Cal B}}}
{\Cal O}(-\lambda)
\,.\tag 5.2$$
Formulas (5.1)-(5.2) yield an isomorphism
$\widetilde{\Cal D}_\lambda^{opp} \iso
\widetilde{\Cal D}_{-\lambda-2\rho}$. 
Taking global sections
we get canonical algebra isomorphisms
$$\tau:\,\widetilde{U}_\lambda^{opp} \simeq
\Gamma({\Cal B},\,\widetilde{\Cal D}_\lambda^{opp}) \iso
\Gamma({\Cal B},\,\widetilde{\Cal D}_{-\lambda-2\rho})
\simeq\widetilde{U}_{-\lambda-2\rho}
\,.\tag 5.3$$

It was shown in [BB3, 3.2.1] that  isomorphisms (5.3) are obtained
by specializing at $\lambda$ a "universal" isomorphism $\widetilde{U}^{opp} \iso
\widetilde{U}$.This way one deduces the following result.

\proclaim{Lemma 5.4} (i) There is a principal anti-involution
$\tau: \widetilde{U} \to \widetilde{U}$ that, for any $\lambda$,
induces the isomorphism $\widetilde{U}_\lambda^{opp}\iso
\widetilde{U}_{-\lambda-2\rho}$ of (5.3).

(ii) The anti-involution
$\tau$ induces anti-involutions
$\tau_{{\frak g}}$ and $\tau_{{\frak h}}$
on the subalgebras $ \Ug \subset\widetilde{U}$ and
${\Bbb C}[{\frak h}^*]\subset\widetilde{U}$, respectively.
The  anti-involutions
$\tau_{{\frak g}}$ and $\tau_{{\frak h}}$ are given,
on generators $x\in {\frak g} \subset 
\Ug$ and $h\in{\frak h} \subset {\Bbb C}[{\frak h}^*]$,
by the formulas:
$$x\mapsto -x\quad,\quad h\mapsto -h-2\cdot\rho(h)\,.$$

(iii) The Harish-Chandra isomorphism $Z({\frak g}) \iso {\Bbb C}[{\frak
h}^*]^W$
intertwines the restriction of $\tau_{{\frak g}}$ to $Z({\frak g})$
with the restriction of $\tau_{{\frak h}}$ to ${\Bbb C}[{\frak
h}^*]^W$.$\quad\square$
\endproclaim

Recall next that the Weyl group $W$ acts on $\widetilde{U}$ and
the action of $w\in W$ induces an isomorphism of completed
algebras $w: \widetilde{U}_{w\cdot \mu} \iso\widetilde{U}_\mu .$
Let $w_0$ be the longest element of $W$. Observe that, for any $\lambda$,
we have
$w_0(\lambda)-2\rho=w_0(\lambda+\rho)-\rho =w_0\cdot\lambda\,.$
Therefore, composing the automorphism of $\widetilde{U}$
induced by $w_0$ with isomorphism  (5.3) we obtain the following isomorphisms
$$\Gamma({\Cal B},\,\widetilde{\Cal D}_{-w_0(\lambda)})\, \iso\,
\widetilde{U}_{-w_0(\lambda)}
\,\,\;{\buildrel {\tau}\over\simeq}\,\,\;
\widetilde{U}_{w_0(\lambda)-2\rho}^{opp}\;
=\;\widetilde{U}_{w_0\cdot\lambda}^{opp}
\,\,\;{\buildrel {w_0}\over\simeq}\,\,\;
\widetilde{U}_\lambda^{opp}\,.\tag 5.5$$

Thus, for any $\widetilde{\Cal D}_{-w_0(\lambda)}$-module
${\Cal N}$, the space  $\Gamma({\Cal B},\,
{\Cal N})$ acquires via (5.5) a natural right 
$\widetilde{U}_\lambda$-module structure. Similarly,
the assignment 
${\Cal M} \mapsto \Gamma({\Cal B} \times {\Cal B},\,
{\Cal M})$   gives a functor from the category of
left  $\widetilde{\Cal D}_\lambda\boxtimes
\widetilde{\Cal D}_{-w_0(\lambda)}$-modules to the category
of $\widetilde{U}_\lambda$-bimodules. Since
both $\lambda$ and $-w_0(\lambda)$  are $\rho$-dominant
weights,
this functor is exact, by 
Theorem 2.2. It will be denoted by $\G$ in the future.

Next, recall that a 
${\Cal D}$-module  ${\Cal M}$  on ${\Cal B} \times {\Cal B}$
is said to be $G$-equivariant if ${\Cal M}$  is a 
$G$-equivariant sheaf of ${\Cal O}_{\Cal{B}\times\Cal{B}}$-modules and the 
differential of the $G$-action coincides with the ${\frak g}$-action on 
${\Cal M}$ arising from the ${\Cal D}$-module structure via the imbedding 
${\frak g} \hookrightarrow$
$\widetilde{\Cal D}_\lambda \boxtimes 
\widetilde{\Cal D}_{-w_0(\lambda)}$. Let  
${\Cal H}{\Cal C}_{\lambda}$ be the category of sheaves of complete 
coherent  $\widetilde{\Cal D}_\lambda \boxtimes 
\widetilde{\Cal D}_{-w_0(\lambda)}$-modules 
on  ${\Cal B} \times {\Cal B}$  that are $G$-equivariant with respect to the 
diagonal $G$-action  on ${\Cal B} \times {\Cal B}$.  
For any ${\Cal M} \in {\Cal H}{\Cal C}_{\lambda}$,
the adjoint 
${\frak g}$-action on the $\widetilde{U}_\lambda$-bimodule 
$\Gamma({\Cal B} \times {\Cal B},\,{\Cal M})$ is 
locally-finite because of the $G$-equivariance of ${\Cal M}$.   
Thus, $\G({\Cal M})$  
becomes an object of $\widetilde{\HC}_{\lambda}$.  Moreover, the functor 
$\G: {\Cal H}{\Cal C}_{\lambda} \to\widetilde{\HC}_{\lambda}
$  is an equivalence of categories provided $\lambda$
is regular (Theorem 2.2).  

Given an algebra $A$, resp. algebras $A$ and $B$, we write
$A\text{-Mod}$ for the category of left $A$-modules,
resp. $\text{Mod-}A$ or
 $A\text{-Mod-}B$
for the category of right $A$-modules, or left $A\boxtimes B^{opp}$-modules (=
$A-B$-bimodules).  We form the bounded  derived categories
$D^b(\widetilde{U}_\lambda\text{-Mod})\,,$
$D^b(\widetilde{U}_\lambda\text{-Mod}$-$\widetilde{U}_\mu)\,$,
and $D^b(\widetilde{\Cal D}_\lambda \boxtimes 
\widetilde{\Cal D}_\mu\text{-Mod})\,.$ Let
$\D\Mod_\lambda(\widetilde{U})$, $\D\widetilde{\HC}_\lambda ,$ 
and $\D{\Cal H}{\Cal C}_{\lambda}$  be the full triangulated subcategories of the 
categories $D^b(\widetilde{U}\text{-Mod})\,$, (resp.
$D^b(\widetilde{U}_\lambda\text{-Mod}$-$\widetilde{U}_\lambda)$, and
$D^b(\widetilde{\Cal D}_\lambda \boxtimes
\widetilde{\Cal D}_{-w_0(\lambda)}\text{-Mod})$ )
formed by the objects whose 
cohomology belong to the categories  
$\Mod_\lambda(\widetilde{U})$, (resp.
$\widetilde{\HC}_{\lambda}$  and  ${\Cal H}{\Cal C}_{\lambda})$.
The functor 
$\G$  has a natural extension to a derived functor  $R\G: 
\D{\Cal H}{\Cal C}_\lambda \to \D\widetilde{\HC}_{\lambda}$.  

\medpagebreak

We now introduce certain convolution functors on our categories.  First, define 
a functor 
$D^b(\widetilde{U}_\lambda\text{-Mod}$-$\widetilde{U}_\mu)\times
D^b(\widetilde{U}_\mu\text{-Mod}$-$\widetilde{U}_\nu)\longrightarrow$
$D^b(\widetilde{U}_\lambda\text{-Mod}$-$\widetilde{U}_\nu)$
 by the formula  $M \star N = M\ 
{\buildrel L\over\otimes}_{\widetilde{U}_\mu}N$. It is easy to
verify that the functor so defined
restricts
to a functor
$$
\star:\;\, \D\widetilde{\HC}_\lambda \times \D\widetilde{\HC}_\lambda \to
\D\widetilde{\HC}_\lambda\;.$$

Observe next that one may regard any
object $A\in \Mod_\lambda(\widetilde{U})$ as an $\widetilde{U}_\lambda$-module on which 
$\widetilde{U}_\lambda$
 acts through a certain quotient $\widetilde{U}_\lambda/{\Cal
J}_\lambda^n\!\cdot\!\widetilde{U}_\lambda$.
Furthermore, given
$M\in \D\widetilde{\HC}_{\lambda}$, one verifies that if
all the cohomology of an object
$A\in \D\Mod_\lambda(\widetilde{U})$ are killed by some power of the
ideal $\J_\lambda$, then a similar vanishing  holds for the cohomology
of $M \star A$. This allows us to 
define a functor
$\D\widetilde{\HC}_\lambda \times \D\Mod_\lambda(\widetilde{U}) \to 
\D\Mod_\lambda(\widetilde{U})$  by the formula  $\enspace
M, A \mapsto M \star A = M\ 
{\buildrel L\over\otimes}_{\widetilde{U}_\lambda}A$. For any   
$M,N,L \in \D\widetilde{\HC}_{\lambda}$,
 and $A \in \D\Mod_\lambda(\widetilde{U})$, 
there are natural functorial isomorphisms (associativity):
$$
(M \star N) \star L \simeq M \star (N \star L)\quad \text{and}\quad 
M \star (N \star A) \simeq (M \star N) \star A .
$$

To define analogous functors on ${\Cal D}$-modules we
introduce the projections 
$p_{ij}: {\Cal B} \times {\Cal B} \times {\Cal B} \to  {\Cal B} \times {\Cal B}$, 
and  $p_i: {\Cal B} \times {\Cal B} \times{\Cal B} \to {\Cal B}$,\ 
$i,j \in \{ 1,2,3\}$, along the factors not named.  Let  $p^*_{ij}\,,\,p^*_i$  stand
for the corresponding inverse image functors in the category of ${\Cal O}$-{\it sheaves}.
 Given a triple of weights 
$(\lambda,\mu,\nu)$, and two objects ${\Cal M} \in D^b(\widetilde{\Cal D}_\lambda 
\boxtimes \widetilde{\Cal D}_\mu\text{-Mod} )$, and ${\Cal N} \in $
$D^b(\widetilde{\Cal D}_{-\mu -2\rho} \boxtimes \widetilde{\Cal D}_\nu\text{-Mod} )$, 
define an object  ${\Cal M} * {\Cal N}$
$ \in D^b(\widetilde{\Cal D}_\lambda 
\boxtimes \widetilde{\Cal D}_\nu\text{-Mod})$  by the formula:
$$
{\Cal M} * {\Cal N} := (Rp_{13})_{{}_\bullet} (p^*_{12} {\Cal M}\ 
{\buildrel L\over\bigotimes}_{p^*_2\widetilde{\Cal D}_\mu} p^*_{23}
{\Cal N})\,,
\tag 5.6
$$
where $(Rp_{13})_{{}_\bullet}$ stands for the sheaf-theoretic
 derived direct image.  It is easy to verify that (5.6)
restricts to a functor
$*:\,\D{\Cal H}{\Cal C}_{\lambda}\times
\D{\Cal H}{\Cal C}_{\lambda}\to$
$\D{\Cal H}{\Cal C}_{\lambda}\,.$

Using the two projections $\text{pr}_i : {\Cal B} \times{\Cal B} \to
{\Cal B}$,
one also defines a convolution-functor  $\,*:\, D^b(\widetilde{\Cal D}_\lambda 
\boxtimes \widetilde{\Cal D}_\mu\text{-Mod} )
 \times D^b(\widetilde{\Cal D}_{-\mu -2\rho}\text{-Mod})
\longrightarrow$
$D^b(\widetilde{\Cal D}_\lambda\text{-Mod} )$  by the formula
$$
{\Cal M} * {\Cal A}\,:\,= (R\text{pr}_1)_{{}_\bullet}({\Cal M}\ 
{\buildrel L\over\otimes}_{{}_{\text{pr}^*_2\widetilde{\Cal D}_\mu}} \text{pr}^*_{2}
{\Cal A})\,.$$
This convolution restricts to a functor
$*:\,\D{\Cal H}{\Cal C}_{\lambda}\times
\D\Mod_\lambda(\widetilde{\Cal D}_\lambda) \longrightarrow$
$\D\Mod_\lambda(\widetilde{\Cal D}_\lambda)$.

The following observation will be useful for us in the next section.

\remark{Remark 5.7} Any object $M\in \D\widetilde{\HC}_{\lambda}$
is completely determined by the corresponding convolution functor
$$M\star\,:\; \D\Mod_\lambda(\widetilde{U}) \to 
\D\Mod_\lambda(\widetilde{U})\quad,\quad
A \mapsto M \star A \,.$$
To see this, we have first to enlarge the category
$\Mod_\lambda(\widetilde{U})$ by adjoining all 
$\widetilde{U}_{\lambda}$-modules which are complete
in the $\J_{\lambda}$-adic topology. Write 
$\widehat{\Mod}_\lambda(\widetilde{U})$ for the latter category,
and $\D\widehat{\Mod}_\lambda(\widetilde{U})$ for the corresponding
derived category enlargement of $\D\Mod_\lambda(\widetilde{U})$.
The functor $M\star (\cdot)$
is continuos in the $\J_{\lambda}$-adic topology, hence extends uniquely to
a well-defined functor $M\star:\,\D\widehat{\Mod}_\lambda(\widetilde{U})
\longrightarrow$
$\D\widehat{\Mod}_\lambda(\widetilde{U})$.

Observe that $\widetilde{U}_\lambda\in
\widehat{\Mod}_\lambda(\widetilde{U})$,
as a left module, and clearly we have
$$M=M\star \widetilde{U}_\lambda = 
{\lim \Sb \\ \longleftarrow\\ n \endSb}
\;\;\Bigl(M\,\star\,
(\widetilde{U}/\J_\lambda^n\!\cdot\!\widetilde{U})\Bigr)\quad,\quad\forall
M\in\D\widetilde{\HC}_{\lambda}\,.
\tag 5.8$$
Notice further
that multiplication by $u\in\widetilde{U}$  on the right gives 
an endomorphism of the
left module $\widetilde{U}/\J_\lambda^n\cdot\widetilde{U}$,
hence induces by functoriality a right $\widetilde{U}_\lambda$-action
on $M\,\star\;
(\widetilde{U}/\J_\lambda^n\cdot\widetilde{U})$. This way,
formula (5.8) recovers $M$ from the functor $M\star(\cdot)$,
as an  $\widetilde{U}_\lambda$-bimodule.
\quad$\square$
\endremark

\medskip

Let  ${\Cal D}(T)$  be the algebra of globally-defined regular differential 
operators on the torus $T$.  View  $S{\frak h}\simeq
{\Bbb C}[{\frak h}^*]$  as the subalgebra 
of  ${\Cal D}(T)$  consisting of $T$-invariant differential operators.  For 
each integer  $n = 1,2,\ldots ,$  define the ${\Cal D}(T)$-module 
${\Cal E}^n = {\Cal D}(T)/{\Cal D}(T)\!\cdot\!\J^n_0$  where  $\J_0$  denotes 
the augmentation ideal in  ${\Bbb C}[{\frak h}^*]$.  Then, ${\Cal E}^n$
gives rise to a sheaf of regular holonomic  ${\Cal D}$-modules on $T$ with 
unipotent monodromy, for $T$ is an affine variety.  These ${\Cal D}$-modules 
form a natural projective system:  ${\Cal E}^1 \twoheadleftarrow 
{\Cal E}^2 \twoheadleftarrow  \ldots\ .$

Now, let ${\Bbb O}$  be the unique open $G$-orbit in ${\Cal B} \times {\Cal B}$, 
$\widetilde{{\Bbb O}} := \pi^{-1}({\Bbb O})$  the inverse image of ${\Bbb O}$ in 
$\widetilde{\Cal B} \times \widetilde{\Cal B}$, and $j: \widetilde{{\Bbb O}} 
\hookrightarrow \widetilde{\Cal B} \times \widetilde{\Cal B}$  the imbedding.
The group $G$ acts {\it freely}    on $\widetilde{{\Bbb O}}$  on the left and the 
space of left cosets, $G\!\setminus\!\widetilde{{\Bbb O}}$,
  is canonically isomorphic to  $T$, due to the
Bruhat decomposition.  Thus, we have a diagram:
$$
\matrix \format \c & \l & \c &\r &\c \\
&&\widetilde{{\Bbb O}} \\
&\mu\swarrow && \searrow \pi \\
\vspace{1\jot}
T=G\!\setminus\!\widetilde{{\Bbb O}}&&&& {\Bbb O} 
\endmatrix \tag 5.9.1
$$

For each $n = 1,2,\ldots ,$ define a
${\Cal D}_{\widetilde{\Cal B} \times\widetilde{\Cal B}}$-module
  of ``multivalued functions" on 
${\widetilde {\Bbb O}} \subset \widetilde{\Cal B} \times \widetilde{\Cal B}$ by
${\Cal O}^n_* = j_*(\mu^* {\Cal E}^n)$.  The sheaf 
$\pi_{{}_\bullet} \widetilde{\Cal O}^n_*$  has a natural $\widetilde{\Cal D} \boxtimes 
\widetilde{\Cal D}$-module structure.  For any weight $\mu \in {\frak h}^*$, 
we set
$$
{\Cal R}_\mu := {\lim \Sb \\ \longleftarrow\\ n \endSb}\;\,
 \Bigl((\widetilde{\Cal D}_\mu 
\boxtimes \widetilde{\Cal D}_{-w_0(\mu)})\,\bigotimes_{\widetilde{\Cal D}
\boxtimes\widetilde{\Cal D}} \pi_{{}_\bullet}{\Cal O}^n_*\Bigr)\, .
\tag 5.9.2$$
We have  ${\Cal R}_\mu \in {\Cal H}{\Cal C}_{\mu}$.  The module  ${\Cal R}_{\mu}$
has the following local description.  
Locally on ${\Bbb O}$, one can 
trivialize the projection 
$\pi: \widetilde{\Cal B} \times \widetilde{\Cal B} \to {\Cal B} \times
{\Cal B}$. Using such a trivialization, one can write
$$
\widetilde{\Cal B} \times \widetilde{\Cal B} \cong {\Cal B} \times {\Cal B} \times T 
\times T\quad  \text{and}\quad  \widetilde{\Cal D}_\mu \boxtimes 
\widetilde{\Cal D}_{-w_0(\mu)} \,\cong\,
 \widetilde{\Cal D}_{\Cal{B}\times\Cal{B}} 
\otimes {\Bbb C}[{\frak h}^*]_\mu
\otimes {\Bbb C}[{\frak h}^*]_{-w_0(\mu)}, 
$$
where  ${\Bbb C}[{\frak h}^*]_\mu$  denotes the $\J_\mu$-adic completion 
of  ${\Bbb C}[{\frak h}^*]$.  Accordingly, locally on ${\Bbb O}$,
one has  ${\Cal R}_\mu = {\Cal O}_{_{\Bbb O}} \otimes {\Bbb C}[{\frak h}^*]_\mu$. 
An element  $1 \otimes a_1 \otimes a_2 \in$
$\widetilde{\Cal
D}_{\Cal{B}\times\Cal{B}}
\otimes {\Bbb C}[{\frak h}^*] \otimes 
{\Bbb C}[{\frak h}^*]$, viewed as an element of  $\widetilde{\Cal D}_\mu 
\boxtimes \widetilde{\Cal D}_{-w_0(\mu)}$ via the local factorization 
above, acts on ${\Cal R}_\mu={\Cal O}_{_{\Bbb O}} \otimes {\Bbb
C}[{\frak h}^*]_\mu$
  as multiplication by $1 \otimes a_1 \cdot 
w_0(\tau(a_2))\,,$ see (5.3) for the meaning of $\tau(a_2)$.

We now define convolution 
functors  $\D{\Cal H}{\Cal C}_\lambda \times \D{\Cal H}{\Cal C}_\lambda 
\to \D{\Cal H}{\Cal C}_{\lambda}$ and $\D{\Cal H}{\Cal C}_\lambda \times 
\D\Mod_\lambda(\widetilde{D}) \to \D\Mod_\lambda(\widetilde{D})$  
by the formula:
$$
{\Cal M} \star {\Cal N} = {\Cal M} * 
{\Cal R}_{w_0(\lambda)-2\rho} * 
{\Cal N}\,,
\tag 5.10 
$$
where  ${\Cal M} \in \D{\Cal H}{\Cal C}_\lambda$ and ${\Cal N}$  is an object of 
either  $\D{\Cal H}{\Cal C}_\lambda$ or $\D\Mod_\lambda(\widetilde{D})$.  Notice
that unlike the $\star$-convolution in (5.10), the $*$-convolution,
  ${\Cal M}*{\Cal N}$, without the middle factor
${\Cal R}_{w_0(\lambda)-2\rho}$ is undefined, for the corresponding
parameters 
$\lambda,\mu,\nu$ in (5.6) do not match in the right way.

It turns out that the convolution (5.10) on ${\Cal D}$-modules corresponds to 
the algebraic $\star$-convolution on $\widetilde{U}_\lambda$-bimodules, i.e., 
we have the following

\proclaim{Proposition 5.11}  
For any ${\Cal M} \in \D{\Cal H}{\Cal C}_\lambda$ 
and ${\Cal N} \in \D{\Cal H}{\Cal C}_\lambda$, resp. ${\Cal N} \in 
\D\Mod_\lambda(\widetilde{U})$,  there is a natural isomorphism
$$
R\G({\Cal M} \star {\Cal N}) \,\simeq\, R\G({\Cal M}) \star 
R\G({\Cal N}),\quad\text{resp.}\enspace
R\Gamma({\Cal M} \star {\Cal N}) \,\simeq\, R\Gamma({\Cal M}) \star 
R\Gamma({\Cal N})\,.
$$
\endproclaim

\demo{Proof}  We have: $\;R\Gamma({\Cal M} \star {\Cal N})=$
$$
= R\Gamma({\Cal M} * 
{\Cal R}_{w_0(\lambda)-2\rho} * {\Cal N}) 
= R\Gamma\bigl({\Cal M} * {\Cal R}_{w_0(\lambda)-2\rho}\bigr)\;
 {\buildrel L\over\bigotimes}_{\widetilde{U}_\lambda}\;
R\Gamma({\Cal N})\,,\tag 5.12
$$
where in the second isomorphism we used that  $\Gamma({\Cal B},\,
\widetilde{\Cal D}_\lambda) 
= \widetilde{U}_\lambda$  and that the isomorphism clearly holds if 
${\Cal M} * {\Cal R}_{w_0(\lambda)-2\rho}$ and ${\Cal N}$  are replaced 
by free $\widetilde{\Cal D}$-modules.
\enddemo

To compute the factor  $R\Gamma({\Cal M} * {\Cal R}_{w_0
(\lambda)-2\rho} * {\Cal N})$ in (5.12)
we apply an extension of the localization theorem
which holds for the derived category of  $\widetilde{\Cal D}_\mu$-modules even if
$\mu$ is {\it not} dominant, see [BB2, \S 12].  Recall first that
we have constructed
an algebra isomorphism 
$\widetilde{U}_\lambda^{opp} \,\simeq\, \widetilde{U}_{-w_0(\lambda)}$,
see (5.3).  It gives rise 
to an exact functor
$\;{\bold {opp}}: \widetilde{U}_{-w_0(\lambda)}$-$\Mod \to$
$\Mod$-$\widetilde{U}_\lambda\;$
that can be extended to derived categories.  On the other hand, the functor 
${\Cal M} \mapsto  {\Cal M} * {\Cal R}_{w_0(\lambda)-2\rho}$  takes 
$D^b(\widetilde{\Cal D}_{-w_0(\lambda)}\text{-Mod} )$ 
to $D^b(\widetilde{\Cal D}_{-\lambda- 2\rho}\text{-Mod} )$.
  Using the canonical isomorphism $\widetilde{\Cal D}_{-\lambda- 2\rho}
\,
\simeq\, \widetilde{\Cal D}_\lambda^{opp}$  we can view the latter functor
as a functor  $*{\Cal R} : 
D^b(\widetilde{\Cal D}_{-w_0(\lambda)}\text{-Mod} )
 \to D^b(\text{Mod-}
\widetilde{\Cal D}_\lambda)$.  Thus, Lemma 5.4 combined
with [BB2, 
theorem of \S 12] yield

{\it For any $\rho$-dominant $\lambda$ the following diagram commutes}:
$$
\matrix \format \c &\quad \c &\quad \c \\
D^b(\widetilde{\Cal D}_{-w_0(\lambda)}\text{-Mod} ) & 
{\buildrel {*{\Cal R}}\over\longrightarrow} 
& D^b(\text{Mod-} \widetilde{\Cal D}_\lambda) \\
\vspace{1\jot}
\big\downarrow R\Gamma && \big\downarrow R\Gamma \\
\vspace{1\jot}
D^b(\widetilde{U}_{-w_0(\lambda)}\text{-Mod} ) & 
{\buildrel {{\bold {opp}}}\over\longrightarrow} 
& D^b(\text{Mod-} \widetilde{U}_\lambda)\;.
\endmatrix
$$
Applying to 
(5.12) a bimodule version of this theorem 
we obtain $R\Gamma({\Cal B}\times {\Cal B},\,
{\Cal M} * {\Cal R}_{w_0(\lambda)-2\rho}) = 
{\bold {opp}}(R\Gamma({\Cal B}\times {\Cal B},\,
{\Cal M}))$.  Hence,  $R\Gamma({\Cal M} \star {\Cal N}) \simeq$
$R\Gamma({\Cal M}) {\buildrel L\over\otimes}_{\widetilde{U}_\lambda}
 R\Gamma({\Cal N})$,
and the proposition follows.\quad $\square$		   
\medpagebreak
 
It should be mentioned, perhaps, that the $*$-convolution
given by formula (5.6) is a $\widetilde{\Cal D}$-analogue of
a more familiar convolution of holonomic 
${\Cal D}_{{\Cal B}\times {\Cal B}}$-modules.
This latter convolution is defined by
$$
{\Cal M} * {\Cal N} := (p_{13})_{{}_*} (p^*_{12} {\Cal M}\ 
{\buildrel L\over\otimes}
p^*_{23}{\Cal N})\;\;\Bigl(=\;\int_{p_{13}}
p^*_{12} {\Cal M}\ 
{\buildrel L\over\bigotimes}_{\Cal O_{{\Cal B}\times {\Cal B}\times {\Cal B}}}
\,p^*_{23}{\Cal N}\Bigr)\,,\tag 5.13
$$
where we are using the same notation
for the 
 projections
$p_{ij}: {\Cal B} \times {\Cal B} \times {\Cal B} \to  {\Cal B} \times
{\Cal B}$,
as in (5.6), and $\int_{p_{13}}$ stands for a 
direct image of ${\Cal D}$-modules. Notice that (5.13)
goes
under the Riemann-Hilbert correspondence to the convolution
on $D^b({\Cal B}\times {\Cal B})$ defined
by  the expression in the middle of (5.13).
\medpagebreak

To proceed further we have to recall some generalities.
Recall that there are Verdier duality functors $\DD$
both on the category of
holonomic ${\Cal D}$-modules and on the category of perverse sheaves,
and that these two functors go to each other under the Riemann-Hilbert
correspondence. It is well known that convolution (5.13)
commutes with Verdier duality. This is no longer  the case for the
convolution (5.6) of $\widetilde{\Cal D}$-modules, due to the fact
that $\widetilde{\Cal B}$ is not compact.  
In more detail, let us restrict to the special case
$\lambda=0$, the only case we need. 
 Verdier duality takes an object of $\D{\Cal H}{\Cal C}_0$
into (in general) a  direct limit of objects of
 $\D{\Cal H}{\Cal C}_0$, an {\it ind}-object
of $\D{\Cal H}{\Cal C}_0$. One can easily
find the commutation relation between convolution (5.6) and
Verdier duality on $\D{\Cal H}{\Cal C}_0$ using the fact that
all objects of this category 
are smooth along the fibers of
the projecion $\widetilde{\Cal B}\times \widetilde{\Cal B}
\to{\Cal B}\times {\Cal B}.$ We have the following result
whose proof is left to the reader.

\proclaim{Lemma 5.14} For any ${\Cal M}\,,\,{\Cal N}\in
\D{\Cal H}{\Cal C}_0$ there is a functorial isomorphism:
$$(\DD{\Cal M}) * (\DD{\Cal N})\;\simeq\;
\DD({\Cal M} * {\Cal N})[\dim_{\Bbb C} T]\;,$$
where $[\dim_{\Bbb C}T]$ stands for the shift in the derived category
by the dimension of the fiber of the projection $\widetilde{\Cal B}
\to{\Cal B}.\quad\square$
\endproclaim

Next, let ${\Cal M} \mapsto {\Cal M}^t$ be the functor on 
$\D{\Cal H}{\Cal C}_0$
induced by the flip of the two factors of the manifold
${\Cal B}\times {\Cal B}$. This functor clearly commutes with Verdier
duality, and we let $\Dt$ denote the composition functor
${\Cal M} \mapsto \DD({\Cal M}^t)=(\DD{\Cal M})^t\,.$
With this understood, one has the following standard result.

\proclaim{Lemma 5.15} For any ${\Cal M}\,,\,{\Cal N}\,,\,{\Cal L}\in
\D{\Cal H}{\Cal C}_0$ there is a functorial isomorphism:
$$
\Hom\bigl({\Cal M} * {\Cal N}, \,{\Cal L}\bigr)\;\simeq\;
\Hom\Bigl({\Cal M} , \,{\Cal L} * (\Dt{\Cal N})\Bigr)\,.\quad\square
$$
\endproclaim
\medpagebreak

In the remainder of this section we are going to formulate a conjecture
relating the various convolution functors introduced above to Koszul
duality. More precisely, we will be dealing with an extension
of the Koszul duality considered in [BGS] to the Harish-Chandra setup,
see [S3].

Write $D^b_G({\Cal B}\times {\Cal B})$ for the $G$-equivariant derived
category on ${\Cal B}\times {\Cal B}$, as defined in [BeLu].
Let $\D_G^{^{\text{mix}}}({\Cal B}\times {\Cal B})$ be its {\it mixed}
version (see [BGS, \S 4.3] where it is referred to as a
{\it graded} version). On the other hand, let 
${\Cal H}{\Cal C}_0^{^{\text{mix}}}$ denote the {\it mixed}
version of the category ${\Cal H}{\Cal C}_0$.
The following
result,
conjectured in [BG],
is an extension of the main theorem of [BGS] to an equivariant setup.

\proclaim{Theorem 5.16} There is a {\bf contravariant} equivalence of triangulated
categories 
$$\K : \D_G^{^{\text{mix}}}({\Cal B}\times {\Cal B})\,\iso\,
\D{\Cal H}{\Cal C}_0^{^{\text{mix}}}\,,$$
that takes pure perverse sheaves in $\D_G^{^{\text{mix}}}({\Cal B}\times {\Cal B})$
to indecomposable projectives in 
${\Cal H}{\Cal C}_0^{^{\text{mix}}}\quad\square$.
\endproclaim

Unfortunately, no complete proof of this theorem has  been so far
written down (although see [S3]). The  remarks below should
hopefully
help
the reader to get a better understanding of the result.

\remark{Remarks 5.17} {\bf (a)} Let $\Per_{_G}({\Cal B}\times {\Cal B})$
denote the abelian category of $G$-equivariant perverse sheaves
on ${\Cal B}\times {\Cal B}$. This is an abelian subcategory
of $D^b_G({\Cal B}\times {\Cal B})$, yet, the natural functor
$D^b\bigl(\Per_{_G}({\Cal B}\times {\Cal B})\bigr)$
$ \,\to\,D^b_G({\Cal B}\times {\Cal B})$ is not an equivalence.

{\bf (b)} The category $\Per_{_G}({\Cal B}\times {\Cal B})$ is known
to be equivalent to the (proper) subcategory of the
category ${\Cal O}_0$ (= category  ${\Cal O}_\lambda$
with $\lambda=0$) formed by  finitely generated $\Ug$-modules
$M$
such that

\noindent
$\bullet\quad$ The augmentation ideal $Z_+({\frak g})=Z({\frak g})\cap
\Ug_+$ annihilates $M$;

\noindent
$\bullet\quad$ $\text{U}\!{\frak b}$-action on $M$ is locally finite;

\noindent
$\bullet\quad$ $\text{U}\!{\frak h}$-action on $M$ is diagonalizable.
\vskip 1pt

{\bf (c)} The category ${\Cal O}_0$, in its turn, is known
to be equivalent to the (proper) subcategory of the
category ${\Cal H}{\Cal C}_0$ formed by the $\widetilde{\Cal D}$-modules
with the trivial monodromy along the fiber of the projection
$\widetilde{\Cal B} \to {\Cal B}$, 
the first factor of the
map $\widetilde{\Cal B}\times\widetilde
{\Cal B} \to {\Cal B}\times {\Cal B}$.$\quad\square$

{\bf (d)} In [BGS] we used instead of the contravariant duality
$\K$, a covariant duality that sends pure perverse sheaves 
in $\D_G^{^{\text{mix}}}({\Cal B}\times {\Cal B})$
to indecomposable injectives in 
${\Cal H}{\Cal C}_0^{^{\text{mix}}}$. The duality $\K$ used here
is obtained from that of [BGS] via composition with Verdier
duality $\DD.\quad\square$.
\endremark
\medpagebreak

Our definitions of the $*$-convolution, resp. $\star$-convolution,
extend verbatim to the category 
$\D_G^{^{\text{mix}}}({\Cal B}\times {\Cal B})$, resp.
$\D{\Cal H}{\Cal C}_0^{^{\text{mix}}}$, and we have

\proclaim{Conjecture 5.18} The Koszul duality functor $\K$ intertwines
the $*$-convolution (5.14) on $\qquad$ $\D_G^{^{\text{mix}}}({\Cal
B}\times {\Cal B})$ with the $\star$-convolution on
$\D{\Cal H}{\Cal C}_0^{^{\text{mix}}}$, that is,
$$\K({\Cal M} * {\Cal N})=\K({\Cal M})\star
\K({\Cal N})\quad,\;\;\forall
{\Cal M}\,,\,{\Cal N}\in\D_G^{^{\text{mix}}}
({\Cal B}\times {\Cal B})\,.$$
\endproclaim
\vskip 1pt

Note that Decomposition theorem [BBD] implies that
$*$-convolution of pure perverse sheaves is pure. On the other hand,
we will see in the next section that  $\star$-convolution
of two projectives in ${\Cal H}{\Cal C}_0$ is again a projective.
This gives a supporting evidence for our conjecture.

We need some more notation. Let $i: \Delta \hookrightarrow 
{\Cal B}\times
 {\Cal B}$
denote the diagonal, and let $\widetilde{\Delta}$ be its inverse image
under the projection $\widetilde{\Cal B}\times \widetilde{\Cal B}\to$
${\Cal B}\times {\Cal B}$. For any point $(x_1,x_2) \in \widetilde{\Delta}$
$\subset \widetilde{\Cal B}\times \widetilde{\Cal B}$
there is a unique element $t= t(x_1,x_2)\in T$ such that the
point $(x_1\cdot t, x_2)$ belongs to the diagonal of
$\widetilde{\Cal B}\times \widetilde{\Cal B}$. The assignment
$(x_1,x_2)\mapsto t(x_1,x_2)$ gives a map $\nu: \widetilde{\Delta}
\to T$, which is analogous to the map $\mu$ in (5.9.1).
Let  $\widetilde{i}: \widetilde{\Delta} \hookrightarrow 
\widetilde{\Cal B}\times \widetilde{\Cal B}$
denote the imbedding. We define the
following objects of the categories 
$\D_G^{^{\text{mix}}}({\Cal B}\times {\Cal B})$ and
${\Cal H}{\Cal C}_0^{^{\text{mix}}}$, respectively:
$$\overline{\Cal L}_\Delta := i_*{\Cal O}_\Delta\quad,\quad
{\Cal L}_\Delta := {\lim \Sb\\ \longleftarrow\\ n\endSb}\;
\Bigl((\widetilde{\Cal D}_0
\boxtimes \widetilde{\Cal D}_0) \bigotimes_{\widetilde{\Cal D} 
\boxtimes \widetilde{\Cal D}} \pi_\bullet (\widetilde{i}_*\nu^*{\Cal E}^n)\Bigr)
\,.
\tag 5.19$$
These objects are the units with respect to $*$-convolutions, that is,
for any $\overline{\Cal M}\in \D_G^{^{\text{mix}}}({\Cal B}\times {\Cal B})$
and
${\Cal M}\in \D{\Cal H}{\Cal C}_0^{^{\text{mix}}}$,
we have
$$\overline{\Cal M} * \overline{\Cal L}_\Delta = \overline{\Cal M} =
\overline{\Cal L}_\Delta * \overline{\Cal M}
\quad,\quad
{\Cal M} * {\Cal L}_\Delta = {\Cal M} =
{\Cal L}_\Delta * {\Cal M}\,.
\tag 5.20
$$

Further, there is the unique open $G$-orbit
$j: {\Bbb O} \hookrightarrow {\Cal B}\times {\Cal B}$. We define the
following objects of $\D_G^{^{\text{mix}}}({\Cal B}\times {\Cal B})$:
$$ \overline{\Cal R}_! := j_!{\Cal O}_{_{\Bbb O}}\quad,
\quad\overline{\Cal R}_*:= j_*{\Cal O}_{_{\Bbb O}}\,.$$
These objects are $*$-inverse to each other, i.e., it is easy to check that
$$\overline{\Cal R}_! * \overline{\Cal R}_* = \overline{\Cal R}_* * \overline{\Cal R}_!
=\overline{\Cal L}_\Delta\,. \tag 5.21
$$
Next, recall  diagram (5.9.1), write $\widetilde{j}_!$ for the direct image with compact
support functor 
corresponding
to the open imbedding 
$\widetilde{j}: \widetilde{\Bbb O} 
\hookrightarrow \widetilde{\Cal B}\times \widetilde{\Cal B}$.
We define
the
following object of the category
${\Cal H}{\Cal C}_\lambda^{^{\text{mix}}}$
$$
{\Cal R}_! := {\lim \Sb\\ \longleftarrow\\ n\endSb}\;
\Bigl((\widetilde{\Cal D}_\lambda
\boxtimes \widetilde{\Cal D}_{-w_0(\lambda)}) \bigotimes_{\widetilde{\Cal D} 
\boxtimes \widetilde{\Cal D}} \pi_\bullet\widetilde{j}_!(\mu^*{\Cal E}^n)\Bigr)\,.
\tag 5.22$$
Then, one can show that
the sheaf ${\Cal R}_! * {\Cal R}_{w_0(\lambda)-2\rho}$  is 
supported on the diagonal in ${\Cal B} \times {\Cal B}$ and, moreover,
the corresponding sheaf on the diagonal is the indecomposable
projective local system with the monodromy representation
(viewed as ${\Bbb C}[{\frak h}]$-module) isomorphic to
${\Bbb C}[{\frak h}]_{-\lambda-2\rho}$. It follows that
the object ${\Cal R}_!$ is the $*$-inverse of
${\Cal R}_{w_0(\lambda)-2\rho}$, i.e., one has
$${\Cal R}_! * {\Cal R}_{w_0(\lambda)-2\rho} = 
{\Cal R}_{w_0(\lambda)-2\rho} * {\Cal R}_!
={\Cal L}_\Delta\,. \tag 5.23
$$

If $\lambda =0$, as we will assume below, we will  write
${\Cal R}_*$ instead of ${\Cal R}_{w_0(\lambda)-2\rho}$,
and view it as an object of 
$\D{\Cal H}{\Cal C}_0^{^{\text{mix}}}$, using a twist
by ${\Cal O}(2\rho)$.
\medskip

\noindent
{\bf {Remark. }} An important  motivation  for Theorem 5.16 comes from
the
equation $\K(\overline{\Cal L}_\Delta) ={\Cal R}_!$, see (5.25) below,
in view of the following. The object
 $\overline{\Cal L}_\Delta$ is a simple $G$-equivariant
perverse sheaf supported on the diagonal of ${\Cal B}\times {\Cal B}$,
and one has:
$$Ext^\bullet_{D^b_G({\Cal B}\times {\Cal B})}(\overline{\Cal L}_\Delta,\, 
\overline{\Cal L}_\Delta) \simeq H^\bullet_G({\Cal B}) \simeq
S{\frak h}^*\,.$$
This corresponds, on the other side of Koszul duality, to 
the fact that ${\Cal R}_!$ is a projective object of 
the category ${\Cal H}{\Cal C}_0^{^{\text{mix}}},$ and moreover
one has an 
isomorphism: $\End_{{\Cal H}{\Cal C}_0^{^{\text{mix}}}}({\Cal R}_!)
\simeq S{\frak h}$.
Similarly, let $\C_{{\Cal B}\times {\Cal B}}$ be
the constant sheaf, viewed as a
 simple object of $D^b_G({\Cal B}\times {\Cal B})$,
and let ${\Cal P}_\Delta$ be the indecomposable projective
cover of ${\Cal L}_\Delta$ in ${\Cal H}{\Cal C}_0^{^{\text{mix}}}$.
Then one has graded algebra isomorphisms:
$$Ext^\bullet_{D^b_G({\Cal B}\times {\Cal B})}(
\C_{{\Cal B}\times {\Cal B}}\,,\,\C_{{\Cal B}\times {\Cal B}}
) \simeq H^\bullet_G({\Cal B}\times {\Cal B})
\simeq
S{\frak h}^*\otimes_{S({\frak h}^*)^W}S{\frak h}^*
\simeq\End_{{\Cal H}{\Cal C}_0^{^{\text{mix}}}}(
{\Cal P}_\Delta\,,\,{\Cal P}_\Delta)\,.\qed
$$

\medpagebreak
\medskip

Assuming Conjecture 5.18 holds true, we
establish the following
commutation relation between the Koszul duality functor $\K$ and
the Verdier duality functor $\DD$.

\proclaim{Theorem 5.24} If Conjecture 5.18 holds then,
for any ${\Cal M}\in
\D_G^{^{\text{mix}}}({\Cal B}\times {\Cal B})$, one has a functorial
isomorphism
$$\K(\Dt{\Cal M})\;\simeq\;
\Dt\Bigl({\Cal R}_* * \K({\Cal M} * \overline{\Cal R}_!\Bigr)
\qquad\text{in}\quad\D{\Cal H}{\Cal C}_0^{^{\text{mix}}}\,.
$$
\endproclaim

\demo{Proof} We begin with the following two easy isomorphisms, which we
leave for the reader to prove
$$ \K(\overline{\Cal L}_\Delta) ={\Cal R}_!\quad,\quad
\K(\overline{\Cal R}_!) ={\Cal L}_\Delta\,.
\tag 5.25$$
Now, for any ${\Cal M}\,,\,{\Cal N}\in
\D_G^{^{\text{mix}}}({\Cal B}\times {\Cal B})$ we have
$$\aligned
&\Hom({\Cal N}\,,\,{\Cal M}) = \quad\text{by (5.20)}\\
&\Hom({\Cal N}\,,\,\overline{\Cal L}_\Delta * {\Cal M}) =  \quad\text{by
Lemma 5.15}\\
&\Hom\bigl({\Cal N} * (\Dt{\Cal M})\,,\,\overline{\Cal L}_\Delta\bigr) = 
\quad\text{by Thm. 5.16}\\
&\Hom\Bigl((\K\overline{\Cal L}_\Delta)\,,\,\K\bigl({\Cal N} * (\Dt{\Cal M})\bigr)\Bigr)
 =  \quad\text{by Conjecture 5.18}\\
&\Hom\Bigl(\K\overline{\Cal L}_\Delta\,,\,\K{\Cal N} * {\Cal R}_* * \K(\Dt{\Cal M})\Bigr)
 =  \quad\text{by Lemma 5.15}\\
&\Hom\Bigl((\K\overline{\Cal L}_\Delta) * (\Dt\cdot\K\cdot\Dt{\Cal M})\,,\,
(\K{\Cal N}) * {\Cal R}_*\Bigr) =  \quad\text{by (5.20)}\\
&\Hom\Bigl((\K\overline{\Cal L}_\Delta) * (\Dt\cdot\K\cdot\Dt{\Cal M})\,,\,
(\K{\Cal N}) * {\Cal R}_* * {\Cal L}_\Delta\Bigr)
 =  \quad\text{by Conjecture 5.18}\\
&\Hom\Bigl((\K^{-1}\K{\Cal N}) * (\K^{-1}{\Cal L}_\Delta)\,,\,
\K^{-1}(\K\overline{\Cal L}_\Delta * (\Dt\cdot\K\cdot\Dt{\Cal M})\bigr)\Bigr)
 =  \\
&\Hom\Bigl({\Cal N} * (\K^{-1}{\Cal L}_\Delta)\,,\,
\K^{-1}\bigl(\K(\overline{\Cal L}_\Delta) * (\Dt\cdot\K\cdot\Dt{\Cal M})\bigr)\Bigr) = 
\quad\text{by (5.25)} \\
&\Hom\Bigl({\Cal N} * \overline{\Cal R}_!\,,\,
\K^{-1}\bigl(\K(\overline{\Cal L}_\Delta) * (\Dt\cdot\K\cdot\Dt{\Cal M})\bigr)\Bigr)
 =  \quad\text{convolve with $\overline{\Cal R}_*$, use  (5.21)}\\
&\Hom\Bigl({\Cal N}\,,\,
\K^{-1}\bigl(\K\overline{\Cal L}_\Delta * (\Dt\cdot\K\cdot\Dt{\Cal M})\bigr)
* \overline{\Cal R}_*\Bigr)\,.
\endaligned
$$
Since these isomorphisms hold for any ${\Cal N}$, we deduce that 
$${\Cal M}\,=\,\K^{-1}\Bigl(\K\overline{\Cal L}_\Delta * (\Dt\cdot\K\cdot\Dt{\Cal M})\Bigr)
* \overline{\Cal R}_*\,.$$
Hence, convolving each side with $\overline{\Cal R}_!$, using (5.21)
and applying $\K$ we get
$$\K({\Cal M} * \overline{\Cal R}_!)\,=\,
(\K\overline{\Cal L}_\Delta) * (\Dt\cdot\K\cdot\Dt{\Cal M})\,.
$$
Thus, using (5.23), (5.25) we calculate
$(\K\overline{\Cal L}_\Delta) * {\Cal R}_* = {\Cal R}_! * {\Cal R}_* = {\Cal L}_\Delta\,.$
Hence, convolving each side with ${\Cal R}_*$, we obtain
$$\Dt({\Cal R}_* * \K({\Cal M} * \overline{\Cal R}_!))\,=\,
\K(\Dt{\Cal M})\,.$$
The theorem is proved $\quad\square$
\enddemo
\bigpagebreak
\subheading{6.\ Projective functors}
\vskip .15in

Fix a $\rho$-dominant integral weight $\lambda$ and write
$\pr_\la$ for the projection of a locally finite
$Z({\frak g})$-module to its $\la$-isotypic component (see after Lemma 1.3).

Following [BeGe, \S3.1], we say that a functor 
$\Phi:
\Mod_\la(\Ug) \to \Mod_\la(\Ug)$  
is a {\it projective functor} if $\Phi$ is a direct summand
of the functor $M\mapsto \pr_\la(E\otimes_{\Bbb C} M)$,
for a certain finite-dimensional ${\frak g}$-module $E$.
Any such functor can be written in terms of $\star$-convolution as
follows. Given $E$ as above, make $E \otimes_{\Bbb C}
\Ug$ 
into a $\Ug$-bimodule by the formulas 
$$
x \cdot (e \otimes u) = x \cdot e \otimes u + e \otimes x \cdot u\quad \text{and}\quad 
(e \otimes u) \cdot x = e \otimes u \cdot x,\quad x \in {\frak g},\ 
e \in E,\ u \in \Ug.$$
Similar formulas make $E \otimes_{\Bbb C}
\Ug_\la$ into a $\Ug$-$\Ug_\la$-bimodule, and for any
$M\in
\Mod_\lambda(\Ug)$, we clearly have
$$\Phi(M)=
\pr_\la(E\otimes_{\Bbb C} M)=\pr_\la\Bigl((E\otimes_{\Bbb
C}\Ug)\,
\otimes_{\Ug}\, M\Bigr)
=\Bigl(
\pr_\la(E\otimes_{\Bbb C}\Ug_{_{\la}})\Bigr)\,\otimes_{{\Ug}_{_{\la}}} M\,.$$
We see that our functor $\Phi$
is a direct summand of the
$\star$-convolution with the bimodule
$\pr_\la(E\otimes_{\Bbb C}\Ug_{_{\la}})\in \HC_{\la}$.
Notice that no (higher) {\it derived} 
 tensor product is required above,
because the bimodule 
$\pr_\la(E\otimes_{\Bbb C}\Ug_{_{\la}})$ 
is a projective right $\Ug_{_{\la}}$-module, as a direct
summand of the free right $\Ug_{_{\la}}$-module
$E\otimes_{\Bbb C}\Ug_{_{\la}}$. 

It is known that the category $\HC_{\la}$ has enough projectives,
 cf. [BGG].
Furthermore, it was shown in [BeGe, \S\S2.2, 4.1] that any  $\Ug_{_{\la}}$-bimodule
of  the form 
$\pr_\la(E\otimes_{\Bbb C}~\Ug_{_{\la}})$ is a projective in the
category
$\HC_{\la}$, and conversely, any projective in $\HC_{\la}$
is a direct summand of $\pr_\la(E\otimes_{\Bbb C}\Ug_{_{\la}})$,
for an appropriate finite dimensional $E$. Convolution with such a
direct summand is by definition a projective functor. Moreover,
we have the following
\proclaim{Proposition 6.1} (i) For any projective functor $\Phi$ there
exists a uniquely determined projective object $P_\Phi\in \HC_{\la}$ such that
$$
\Phi(M) = P_\Phi \star M\quad\text{where}\quad
P_\Phi ={\lim \Sb \\ \longleftarrow\\ n \endSb}\;\,
\Phi\bigl(\Ug/\I_\la^n\!\cdot\!\Ug\bigr)\,;
$$
The assignment 
$\Phi \mapsto \PP_\Phi$ sets up a (1-1)-correspondence 
between projective functors and 
projectives in  $\HC_{\la}$.  

(ii) Any projective functor is exact.

(iii) Composition of projective functors is a projective functor.
\endproclaim

\demo{Proof} Part (i) follows from an analogue of Remark 5.7
(for $\Ug$ instead of $\widetilde{U}$).  Part (ii) follows from
the exactness of the functor 
$\pr_\la(E\otimes_{{\Bbb C}} \bullet)$. Part (iii) is clear from
definition.\quad$\square$
\enddemo

For any projective functor $\Phi$, there is a projective functor  $\Phi^\dag$  which 
is both the left and the right adjoint of  $\Phi$, see [BeGe, Lemma 3.2(v)].  

Next we extend projective functors to the category  $\Mod_\lambda(\widetilde U)$ 
as follows.  To any ${\Ug_\la}$-bimodule $M \in \HC_{\la}$
associate the $\widetilde{U}_\lambda$-bimodule 
$$\widetilde{M} \,:=\,
\widetilde{U}_\lambda \otimes_{\widetilde{U}_\lambda^{W_\lambda}} M 
\otimes_{\widetilde{U}_\lambda^{W_\lambda}} \widetilde{U}_\lambda \in 
\widetilde{\HC}_\lambda\,.$$
see Lemma 1.3.  If $M$ is a projective in 
$\HC_{\la}$, then $\widetilde M$  is a projective in 
$\widetilde{\HC}_\lambda$ (by adjunction of $\Res_\lambda$ and $\widetilde{U}_\lambda
\otimes (\cdot))$.  Given a projective functor  $\Phi$ on $\Mod_{\la}(\Ug)$, 
define a functor  $\widetilde{\Phi}$  on $\Mod_\lambda(\widetilde{U})$ by
$$
\widetilde{\Phi}(M) \,:= \,
{\buildrel\widetilde{\qquad\quad}\over {\overbrace {P_\Phi\otimes_{\Ug_{\la}}\Res_\lambda M}}}
= \widetilde{\PP}_\Phi \otimes_{\widetilde{U}_\lambda} M.
$$
Clearly, $\widetilde{\Phi}:\Mod_\lambda(\widetilde{U})
\to \Mod_\lambda(\widetilde{U})$  is an exact functor and, moreover, 
$\widetilde{\Phi}(M) = \widetilde{\PP}_\Phi \star M$.  

The following result describes the relation between projective functors and 
${\Cal D}$-modules.

\proclaim{Proposition 6.2} (i)  For any projective  ${\Cal P} \in 
{\Cal H}{\Cal C}_\lambda$  the assignment  ${\Cal M} \mapsto  {\Cal P} \star 
{\Cal M}$  gives an exact functor  {$\Mod_\lambda(\widetilde{\Cal D}) \to 
\Mod_\lambda(\widetilde{\Cal D})$ and there is a projective functor
$\widetilde{\Phi} = \widetilde{\Phi}_{\Cal P}$ on $\Mod_\lambda(\widetilde{U})$
such that the following diagram commutes:
$$
\matrix \format \c &\quad \c &\quad \c \\
\Mod_\lambda(\widetilde{\Cal D}) & {\buildrel {\Gamma}\over\longrightarrow} 
& \Mod_\lambda(\widetilde{U}) \\
\vspace{1\jot}
{\Cal P}\star(\cdot)\big\downarrow {\phantom{PPP}} && \big\downarrow \widetilde{\Phi} \\
\vspace{1\jot}
\Mod_\lambda(\widetilde{\Cal D}) & {\buildrel {\Gamma}\over\longrightarrow} 
& \Mod_\lambda(\widetilde U) \quad .
\endmatrix 
$$

\item{(ii)} For any projective functor $\widetilde{\Phi}$  on $\Mod_\lambda(\widetilde U)$
there exists a projective ${\Cal P} \in {\Cal H}{\Cal C}_\lambda$ such that
$\widetilde{\Phi} = \widetilde{\Phi}_{\Cal P}$, i.e., such that the
above diagram commutes.
\endproclaim

\demo{Proof}  By Proposition 5.11 we have  $\Gamma({\Cal P} \star 
{\Cal M}) = \G({\Cal P}) \star \Gamma(M)$.  If $\lambda$ is 
regular, then $\G$  gives an equivalence of the categories  
${\Cal H}{\Cal C}_\lambda$  and $\widetilde{\HC}_\lambda$. This proves  
the Proposition if $\lambda$ is regular.

To complete the proof of part (i) in the general case choose a {\it regular} 
$\rho$-dominant weight $\mu$.   Using the geometric translation functor 
$\Theta^\mu_\lambda$  we obtain the following equivalences of categories:
$$
{\Cal H}{\Cal C}_\lambda \simeq {\Cal H}{\Cal C}_\mu \simeq \widetilde{\HC}_\mu 
\simeq \HC_{\ma}\quad,\quad \Mod_\lambda(\widetilde{\Cal D}) 
\simeq \Mod_\mu(\widetilde{\Cal D}) \simeq \Mod_\mu(\widetilde{U}) 
\simeq \Mod_{\ma}(\Ug).
$$

The Proposition being already known for regular $\mu$'s, it follows that the 
functor  ${\Cal P} \star (\cdot)$  
is an exact functor on $\Mod_\lambda(\widetilde{\Cal D})$
and that there exists a projective  ${\Cal P}^\dag \in {\Cal H}{\Cal C}_\lambda$
such that the functor  ${\Cal P}^\dag \star(\cdot)$  is both the left and the 
right adjoint of  ${\Cal P} \star(\cdot)$.  Now, the category  
$\Mod_\lambda(\widetilde{U})$  may be viewed as a quotient of the category 
$\Mod_\lambda(\widetilde{\Cal D})$.  By [S1, Lemma 6, p.432], the pair of
adjoint functors  $({\Cal P} \star(\cdot),$ 
${\Cal P}^\dag \star(\cdot))$ descends 
to the pair  $(\G({\Cal P}) \star(\cdot),\,
 \G({\Cal P}^\dag) \star(\cdot))$
  of adjoint (exact) functors on  $\Mod_\lambda(\widetilde U)$.
We have in particular 
$$\aligned
\Hom_{\Mod_\lambda(\widetilde U)}\Bigl(\widetilde{U}_\lambda, 
\,\G({\Cal P}^\dag)\star M\Bigr) =& 
\Hom_{\Mod_\lambda(\widetilde U)}\Bigl(\G({\Cal P}) 
\star \widetilde{U}_\lambda, \,M\Bigr)\\
& \quad= 
\Hom_{\Mod_\lambda(\widetilde U)}\bigl(\G({\Cal P}), \,M\bigr).
\endaligned \tag 6.3
$$
One checks from the construction that isomorphisms
(6.3) still hold if $M$ is taken to be 
a bimodule from  $\widetilde{\HC}_\lambda$  and  $\Hom$'s are taken in 
$\widetilde{\HC}_\lambda$.  But then the functor
$\Hom_{\Mod_\lambda(\widetilde U)}\Bigl(\widetilde{U}_\lambda, 
\,\G({\Cal P}^\dag)\star (\bullet)\Bigr)$ on the left-hand side of 
(6.3) is an exact functor.  Hence, the functor
$\Hom_{\Mod_\lambda(\widetilde U)}(\G({\Cal P}),
\,\bullet)$
on the right-hand 
side is an exact functor again.  Thus, $\G({\Cal P})$  is a projective in
$\widetilde{\HC}_\lambda$  and part (i) follows.

To prove part (ii) it suffice  to show, by part (i), that any projective 
$P \in \widetilde{\HC}_\lambda$ can be written in the form  $P = 
\G({\Cal P})$,  where ${\Cal P}$  is a projective in 
${\Cal H}{\Cal C}_\lambda$.  We may put  ${\Cal P} := \Delta_\lambda P$
where  $\Delta_\lambda$  is the localization functor.  Then ${\Cal P}$  is 
projective, by adjunction, and the result follows from Theorem 2.2(i).
\quad $\square$
\enddemo

As a simple application of our analysis we get

\proclaim{Corollary 6.4} There is a natural isomorphism
 $\G({\Cal R}_!) \simeq 
\widetilde{U}_\lambda$.
\endproclaim
\demo{Proof}  Observe that formula
(5.23) yields
$$
{\Cal R}_! * {\Cal R}_{w_0(\lambda)-2\rho} * {\Cal M} =
{\Cal M}\;,
\quad
\text{for\ any}\ \ {\Cal M} \in \Mod_\lambda(\widetilde{\Cal D}).
$$
This equation  shows that $ \star$-convolution with
${\Cal R}_!$ gives the identity 
functor on $\widetilde{{\Cal H}{\Cal C}}_\lambda$.
On the other 
hand, view the bimodule  $\widetilde{U}_\lambda$  as an object of 
$\widetilde{\HC}_\lambda$. We know that the functor
$\widetilde{U}_\lambda \star (\cdot )$ is the identity 
functor on $\widetilde{\HC}_\lambda$, cf. (5.8). Thus, Proposition
5.11 yields the result.$\quad\square$
\enddemo

Now let $\lambda, \mu$  be integral $\rho$-dominant weights such that 
$W_\lambda \subset W_\mu$  and let $\widetilde{\theta}^-: \Mod_\lambda(\widetilde{U})
\to \Mod_\mu(\widetilde{U})$ be the translation functor.  Repeating the 
definitions, one gets a similar functor  $\widetilde{\theta}^-: 
\widetilde{\HC}_\lambda \to \widetilde{\HC}_\mu$  on bimodules.

\proclaim{Corollary 6.5}  Let $P$ be a projective in $\widetilde{\HC}_\lambda$. 
Then  $\widetilde{\theta}^-P$  is a projective in $\widetilde{\HC}_\mu$  and the following 
diagram of functors commutes:
$$
\matrix \format \c &\quad \c &\quad \c \\
\Mod_\lambda(\widetilde{U}) & {\buildrel {P\star(\cdot)}\over\longrightarrow} 
& \Mod_\lambda(\widetilde{U}) \\
\vspace{1\jot}
\widetilde{\theta}^-\big\downarrow && \big\downarrow \widetilde{\theta}^- \\
\vspace{1\jot}
\Mod_\mu(\widetilde{U}) & {\buildrel {(\widetilde{\theta}^-P) \star 
(\cdot)}\over\longrightarrow} & \Mod_\mu(\widetilde{U}) .
\endmatrix
$$
\endproclaim

\demo{Proof}  Let $\Theta^-: {\Cal H}{\Cal C}_\lambda \longrightarrow 
{\Cal H}{\Cal C}_\mu$  be the geometric translation functor.
By the proof of part (ii) of Proposition 6.2, there is a projective
${\Cal P} \in {\Cal H}{\Cal C}_\lambda$  such that  $P = \G({\Cal P})$.
Then, 
by Proposition 2.8, we have  $\theta^-P = \G(\Theta^-{\Cal P})$.
The proof of part (i) of Proposition 6.2 shows now that 
$\G(\Theta^-{\Cal P})$  is a projective in $\widetilde{\HC}_\mu$.  
Furthermore, one checks easily that the following ${\Cal D}$-module counterpart 
of the above diagram commutes:
$$
\matrix \format \c &\quad \c &\quad \c \\
\Mod_\lambda(\widetilde{{\Cal D}}) & {\buildrel {{\Cal P}\star(\cdot)}\over\longrightarrow} 
& \Mod_\lambda(\widetilde{{\Cal D}}) \\
\vspace{1\jot}
\widetilde{\Theta}^-\big\downarrow && \big\downarrow \widetilde{\Theta}^- \\
\vspace{1\jot}
\Mod_\mu(\widetilde{{\Cal D}}) & {\buildrel {(\widetilde{\Theta}^-{\Cal P}) \star 
(\cdot)}\over\longrightarrow} & \Mod_\mu(\widetilde{{\Cal D}}) .
\endmatrix
$$
The result now follows from Proposition 6.2.\quad $\square$
\enddemo

\medpagebreak

Recall the functor  $\widetilde{\theta}^+_\ell:	 \Mod_\mu(\widetilde{U}) 
\to \Mod_\lambda(\widetilde{U})$, the left adjoint of $\widetilde{\theta}^-$.
We have:

\proclaim{Proposition 6.6}  Assume that for some $M \in \Mod_\lambda(\widetilde{U})$
the adjunction morphism  $\widetilde{\theta}^+_\ell \cdot \widetilde{\theta}^-M \to M$ 
is an isomorphism.  Then, for any projective functor  $\widetilde{\Phi}$ on 
$\Mod_\lambda(\widetilde{U})$, the adjunction morphism 
$\widetilde{\theta}^+_\ell \cdot \widetilde{\theta}^-(\widetilde{\Phi}M) \to 
\widetilde{\Phi} M$  is an isomorphism again.
\endproclaim

\demo{Proof}  Let  $\widetilde{\Phi}^\dag$  be the adjoint of a projective functor 
$\widetilde{\Phi}$.  The pair $(\widetilde{\Phi}, \widetilde{\Phi}^\dag)$  descends, 
by Corollary 6.5 and [S1, Lemma 6, p.432] to an adjoint pair of exact functors 
$(\widetilde{\Psi}, \widetilde{\Psi}^\dag)$ on $\Mod_\mu(\widetilde{U})$.  For 
any  $M,N \in \Mod_\lambda(\widetilde{U})$  we have, by adjunction, 
$$
\aligned
\Hom(\widetilde{\theta}^+_\ell \cdot \widetilde{\theta}^-\widetilde{\Phi} M,\,N) 
&= \Hom(\widetilde{\theta}^-\widetilde{\Phi} M,\,
\widetilde{\theta}^-N)\\
& =
\Hom(\widetilde{\Psi}\widetilde{\theta}^-M,\, \widetilde{\theta}^-N) \\
& = \Hom(\widetilde{\theta}^-M,\, 
\widetilde{\Psi}^\dag\widetilde{\theta}^-N)\\
& = 
\Hom(\widetilde{\theta}^-M,\, \widetilde{\theta}^-\widetilde{\Phi}^\dag
N)\\
& = 
\Hom(\widetilde{\theta}^+_\ell \cdot \widetilde{\theta}^-M,\, \widetilde{\Phi}^\dag N)\,.
\endaligned
$$
\noindent
If  $\widetilde{\theta}^+_\ell \cdot \widetilde{\theta}^-M = M$, then the last $Hom$
can be rewritten as 
$$
\Hom(M, \widetilde{\Phi}^\dag N) = \Hom(\widetilde{\Phi} M, N)\,,
$$
and the  proposition follows. \quad $\square$
\enddemo
\medpagebreak

Next let $\lambda=0$. Since $0$ is a $\rho$-dominant regular
weight, we may identify the categories ${\Cal H}{\Cal C}_0$
and $\widetilde{HC}_0$. There is a standard duality functor
on the category $\widetilde{HC}_0$.  Any projective
functor is known to commute with that duality,
since so does tensoring with a finite
dimensional
representation,  see [BeGe]. 
Moreover, it is known
 that
the standard duality on $\widetilde{HC}_0$ goes into the Verdier
duality on ${\Cal H}{\Cal C}_0$ (this is false for the Harish-Chandra
category over the general {\it real} reductive group: the two
dualities
may, in general,
act differently already on the simple objects). The two dualities
coincide however in the special case of
a {\it complex} reductive group, the case we are interested in.
Then the claim 
 boils down to a similar
result for the category ${\Cal O}$. In the category ${\Cal O}$ case
the result is known,
 although we could not find any written account of it.

Here is a sketch of proof. First, it is immediate to verify
that the two
dualities agree on every simple object of ${\Cal O}$.
Second, each duality is an exact functor, hence,
agreement on simple objects implies that any 
object $M\in {\Cal O}$ is sent by both dualities to
isomorphic objects, call it $M^\dag$.
The remaining (most delicate) part of the proof is to verify that,
for any $M, N\in {\Cal O}$,
the two dualities induce the same maps:
$\Hom_{{\Cal O}}(M, N) \to \Hom_{{\Cal O}}(N^\dag, M^\dag)$. This is equivalent,
by a standard homological algebra, to the claim that,
for any {\it simple} objects $M, N\in {\Cal O}$ and any $i\ge 0$,
the two dualities induce the same maps:
$Ext^i_{{\Cal O}}(M, N) \to Ext^i_{{\Cal O}}(N^\dag, M^\dag)$.
The latter claim is obvious for $i=0$, and can be verified by hand
for $i=1$, using that for $i=1$ one can take $M$ to be a Verma module
instead of a simple module. The general case $i\ge 1$
now follows from
the main theorem of [BGS] which implies that the $Ext$-algebra of
all the simple objects is Koszul, in particular, is generated
by $Ext^0$'s and $Ext^1$'s.\qed\medskip

We now return to our convolutions and note that for
$\lambda=0$
we may assume, twisting by the canonical bundle
if necessary, that both $*$-convolution, and the $\star$-convolution
take the category  $\D{\Cal H}{\Cal C}_0$ into itself and, moreover,
that the object ${\Cal R}$ involved in formula (5.10) is also
an object of ${\Cal H}{\Cal C}_0$. Recall further that the
 Verdier duality functor, $\DD$, is a contravariant exact functor
on ${\Cal H}{\Cal C}_0$. Hence, it takes projective objects into
injective
objects.

\proclaim{Proposition 6.7} For any projective
${\Cal P}\in{\Cal H}{\Cal C}_0$, in the notation of (5.9.1)-(5.9.2) we have 
$$\DD{\Cal P} = {\Cal P} \star {\Cal R}^\vee\quad\text{where}\quad
{\Cal R}^\vee := {\lim \Sb \\ \longrightarrow\\ n \endSb}\;\,
 \Bigl((\widetilde{\Cal D}_0
\boxtimes \widetilde{\Cal D}_0)\,\bigotimes_{\widetilde{\Cal D}
\boxtimes\widetilde{\Cal D}}\,
\pi_{{}_\bullet}j_*\mu^*(\DD{\Cal E}^n)\Bigr)\,.
$$
\endproclaim
\demo{Proof} Let $\Phi_{\Cal P} : M\mapsto {\Cal P}\star M$ be the
projective functor
corresponding to ${\Cal P}$.  It follows from the discussion
of two dualities above the proposition
that the functor $\Phi_{\Cal P}$ on ${\Cal H}{\Cal C}_0$
commutes with Verdier
duality. Therefore, setting $r=\dim T$
and using Lemma 5.14, for any ${\Cal M}\in {\Cal H}{\Cal C}_0$, we find
$$
\aligned
{\Cal P} * {\Cal R} * (\DD {\Cal M}) & =
{\Cal P}\star (\DD {\Cal M})\\
& =
\Phi_{\Cal P}(\DD {\Cal M})\qquad\text{since 
$\Phi$ commutes with Verdier duality}\\
& = \DD\cdot\Phi_{\Cal P}({\Cal M})=
\DD({\Cal P}\star {\Cal M})\\
&=\DD({\Cal P} * {\Cal R} * {\Cal M})
\quad\text{by Lemma 5.14}\\
&=\DD({\Cal P} * {\Cal R}) * \DD({\Cal M})[r]\,.
\endaligned
$$
The above isomorphisms imply, by Remark 5.7, that
$${\Cal P} * {\Cal R} = \DD({\Cal P} * {\Cal R})[r] =
(\DD {\Cal P}) * (\DD {\Cal R})\,,\tag 6.8$$
where in the last equality we have used Lemma 5.14 once more.
Now, it is easy to show that the object 
${\Cal R}^\vee$ is $*$-inverse to 
$\DD{\Cal R}$. Therefore,
convolving both sides of (6.8) with ${\Cal R}^\vee$ yields
$$\DD{\Cal P} = {\Cal P} * {\Cal R} * {\Cal R}^\vee =
{\Cal P} \star {\Cal R}^\vee\,,$$
and the claim follows. $\quad\square$
\enddemo
\medpagebreak

{\bf Tilting $\pmb{\Cal D}$-modules.}
We continue to assume that $\lambda =0$, and use the notation of
\S4. Let ${\frak b}$ be the
 fixed Borel subalgebra, and 
$B\subset G$ the corresponding Borel group.
Since $0$  is regular we may (and will)
identify 
$\Ug$-modules with $\widetilde{\Cal D}$-modules via the equivalences: 
$\Mod_{|0|}(\Ug) \simeq \Mod_0(\widetilde{U}) \simeq
\Mod_0(\widetilde{\Cal D})$.
 Further, there is a natural equivalence, see e.g. [BeGe, \S5],
between the category
of $\text{U}\!{\frak b}$-locally finite $\Ug_{_{|0|}}$-modules,
resp. $B$-monodromic  $\widetilde{\Cal D}$-modules on
${\Cal B}$, see [BB3], and the
category ${\Cal H}{\Cal C}_0$. We will write
$\dd(M)$ for the object of ${\Cal H}{\Cal C}_0$
corresponding to $M\in \Mod_{_{|0|}}(\Ug)$.
As we mentioned in Remark 5.17(c), the category ${\Cal O}_0$
goes under the above equivalence to the subcategory of 
 ${\Cal H}{\Cal C}_0$ formed by the $\widetilde{\Cal D}$-modules
with the trivial monodromy along the fiber of the projection
$\widetilde{\Cal B} \to {\Cal B}$, 
the first factor of the
map $\widetilde{\Cal B}\times\widetilde
{\Cal B} \to {\Cal B}\times {\Cal B}$. 
In particular, for each  $w \in W$, there are objects
$\dd(M_{w\cdot 0})\,,\, \dd(P_{w\cdot 0})\in {\Cal H}{\Cal C}_0$
corresponding to the Verma module and its indecomposable projective
cover, respectively.

Given $\mu \in {\frak h}^*$, let $\text{S}{\frak h}_\mu$
denote the
completion of the Symmetric algebra at the corresponding point.
We introduce the "universal" Verma
module ${\Bbb M} := \Ug/ \Ug\!\cdot\!{\frak n}$,
where ${\frak n} :=[{\frak b} , {\frak b}]$ is the nilradical of 
${\frak b}$.
This module has a natural {\it right}  $\text{S}{\frak h}$-action,
and we set ${\Bbb M}_\mu := {\Bbb M}
\bigotimes_{_{\text{S}\!{\frak h}}}\text{S}{\frak h}_\mu$, the
completeion of ${\Bbb M}$
at the point $\mu$. In particular, for each  $w \in W$,
there is an $\Ug$-module ${\Bbb M}_{w\cdot 0}$,
which is not an object of the category ${\Cal O}_0$. We set
${\Cal M}_{w}:= \dd({\Bbb M}_{w\cdot 0})$ 
denote the corresponding object in
${\Cal H}{\Cal C}_0$, and let ${\Cal P}_{w}$ be
its indecomposable projective
cover in ${\Cal H}{\Cal C}_0$. We note that,
if $w=w_0$, then ${\Cal M}_{w_0} = {\Cal R}_!$ is the Harish-Chandra
module defined in (5.22). In general, ${\Cal M}_{w}$ has a
 pro-unipotent monodromy (i.e., free $\text{S}{\frak h}_{w\cdot 0}$-action)
along the  fiber of the projection
$\widetilde{\Cal B} \to {\Cal B}$ 
(the first factor of the
map $\widetilde{\Cal B}\times\widetilde
{\Cal B} \to {\Cal B}\times {\Cal B}$), and 
$\dd({ M}_{w\cdot 0})$
is obtained from ${\Cal M}_{w}$ by 
"killing" (= taking co-invariants of) the monodromy along the first
factor
(very much the same  way  as  the Verma module
$M_\mu$ is obtained from ${\Bbb M}_\mu$).
Similarly, the projective ${\Cal P}_{w}$
has a flag formed by the various  ${\Cal M}_{y}$'s,
hence has free pro-unipotent monodromy; again,
$\dd(P_{w\cdot 0})$ is obtained from ${\Cal P}_{w}$ by 
"killing" that monodromy. In particular, one has an equality of
multiplicities: 
$$\,[{\Cal P}_{w}\, : \,
{\Cal M}_{y}] =
[\dd(P_{w\cdot 0}) \,:\, \dd(M_{y\cdot 0})]\,.\tag 6.9$$

We will also need Harish-Chandra modules ${\Cal M}_w^\vee$
which are in a sense dual analogues of ${\Cal M}_w$. We first
define the corresponding $\Ug$-modules as follows. For $w \in W$, and
each 
$n=1,2,\ldots,$ we have a well-defined object 
${\Bbb M}_\mu/\J_\mu^n \in {\Cal O}_0$, where
$\J_\mu\subset \text{S}{\frak h}_\mu$ is  the
maximal ideal. Let $\DD({\Bbb M}_\mu/\J_\mu^n)$ be its dual
in ${\Cal O}_0$. There are natural projections (not imbeddings !):
$\DD({\Bbb M}_\mu/\J_\mu^{n+1}) \twoheadrightarrow
\DD({\Bbb M}_\mu/\J_\mu^n)\,,\,n=1,2,\ldots,$ and we set
$\displaystyle{\Bbb M}_\mu^\vee := {\lim \Sb\\ \longleftarrow\\ n\endSb}\;
\DD({\Bbb M}_\mu/\J_\mu^n)$.
Let ${\Cal M}_w^\vee := \dd({\Bbb M}_{w\cdot 0}^\vee)$ be the
corresponding Harish-Chandra modules.
We note that for $w=w_0$ we have:
${\Cal M}_{w_0}^\vee = {\Cal R}_{-2\rho}$ is the object
introduced in (5.9.2), see also (5.23).
\smallskip

It was shown in [BeGe] that, for each  $w \in W$, there
is a unique indecomposable projective functor  $\Phi_w$, see [BeGe,
Thm. 3.3], such that:
$\Phi_w(M_0) = P_{w\cdot 0}$. Hence, for each  $w \in W$,
there is, by Propositions 6.1-6.2, a unique projective 
${\Cal P}_w \in {\Cal H}{\Cal C}_0$,
such that the functor $\Phi_w$ corresponds on
$\widetilde{\Cal D}$-modules to the 
convolution functor  ${\Cal M} \mapsto
{\Cal P}_w * {\Cal R} * {\Cal M}$, where ${\Cal R} := {\Cal
R}_{-2\rho}$
is as above.
 It is not difficult to see
that ${\Cal P}_w$ is the indecomposable
projective cover of ${\Cal M}_w$ (considered 2 paragraphs above),
equivalently, the projective cover
of the simple $\widetilde{\Cal D}$-module
supported on the $G$-diagonal orbit in ${\Cal B}\times {\Cal B}$,
corresponding to the element $w$.

Recall that an object  of the category ${\Cal O}$ is
called a {\it tilting} module, if it is self-dual and
has a Verma-flag. 
Observe that projective functors take
self-dual objects into self-dual ones, and modules 
with Verma-flag into modules 
with Verma-flag, cf. [BeGe]. 
Since $M_{-2\rho}$ is a simple, hence, self-dual Verma
module, it follows that, for each  $w \in W$,
the module
$Q_{w\cdot 0} := \Phi_w(M_{-2\rho})$ is tilting. In fact,
the $\{Q_{w\cdot 0}\,,\,w\in W\},$ are exactly all the
 indecomposable tilting modules in 
the category ${\Cal O}_0$. 

\proclaim{Definition} A finite direct sum of objects
of $\D{\Cal H}{\Cal C}_0$ of the form
${\Cal Q}_w
:= {\Cal P}_w * {\Cal R}\,,\,w \in W$, cf. (5.10),
will be called a {\it tilting} Harish-Chandra module.
\endproclaim

Note that this definition gives a natural mixed structure
on ${\Cal Q}_w$, cf. [BGS, \S4].

\proclaim{Theorem 6.10}  (i) 
For any $w \in W$ , the complex ${\Cal Q}_w := {\Cal P}_w * {\Cal R}$
is actually a  $\widetilde{\Cal D}$-module
(not just an object of the derived category), 
${\Cal Q}_w \in {\Cal H}{\Cal C}_0$.

(ii) For any $y, w \in W$ we have: ${\Cal Q}_y * {\Cal Q}_w
\in {\Cal H}{\Cal C}_0$, is again a tilting module.

(iii) For any $w, y \in W$ the following multiplicity
formula holds: $[Q_w : {\Cal M}_{y}^\vee]=
[P_w : M_{yw_0}]$. Moreover:
$\Hom_{{\Cal H}{\Cal C}_0}({\Cal Q}_w, {\Cal Q}_y)
\simeq \Hom_{{\Cal H}{\Cal C}_0}({\Cal P}_w, {\Cal P}_y)$ and
$\Hom_{{\Cal O}}(Q_{w\cdot 0}, Q_{y\cdot 0})
\simeq \Hom_{{\Cal O}}(P_{w\cdot 0}, P_{y\cdot 0})$.
\endproclaim

Before going into proof we note 
that Proposition 6.2 yields: ${\Cal Q}_w * 
\dd(M_{-2\rho}) = {\Cal P}_w * {\Cal R} *
\dd(M_{-2\rho})$
$=
{\Cal P}_w \star
\dd(M_{-2\rho}) = \dd\bigl(\Phi_w(M_{-2\rho})\bigr)=\dd(Q_{w\cdot 0})$. 
Furthermore, one shows as above 
 that the  Harish-Chandra module $\dd(Q_{w\cdot 0})$
is obtained
from the tilting Harish-Chandra modules ${\Cal Q}_w$ by 
taking coinvariants of the monodromy 
along the  fiber of the projection
$\widetilde{\Cal B} \to {\Cal B}$
(the first factor of the
map $\widetilde{\Cal B}\times\widetilde
{\Cal B} \to {\Cal B}\times {\Cal B}$).

\demo{Sketch of proof of Theorem 6.10}
To prove (i) we exploit an  interpretation of the functor:
$\,{\Cal M} \mapsto {\Cal P}_y * {\Cal R} * {\Cal M}$ as of a projective
functor. Such functors take $\widetilde{\Cal D}$-modules into
$\widetilde{\Cal D}$-modules,
because of the global interpretation via bi-modules over the extended
enveloping
algebra $\widetilde{U}$. Hence, using the notation of (5.19)
we get: ${\Cal Q}_y =
{\Cal Q}_y * {\Cal L}_\Delta = {\Cal P}_w * {\Cal R} * {\Cal L}_\Delta
={\Cal P}_y \star {\Cal L}_\Delta
\in {\Cal H}{\Cal C}_0$, and part (i) is proved.
Further, for any $M \in \Mod_{_{|0|}}(\Ug)$,
one has:
$({\Cal Q}_y * {\Cal Q}_w) * M= {\Cal P}_y * {\Cal R} *{\Cal P}_w *
{\Cal R} * M$
$= ({\Cal P}_y \star {\Cal P}_w) \star M\,$.
Part (ii) now follows from Proposition 6.1(iii) and Remark 5.7.

Observe next that, for any $y \in W$, we
have $\ell(y^{-1}) + \ell(yw_0) = \ell(w_0)$.
It follows that: ${\Cal M}_{y^{-1}}^\vee * {\Cal M}_{yw_0}^\vee
={\Cal M}_{w_0}^\vee = {\Cal R}$. Also, one knows that:
${\Cal M}_y * {\Cal M}_{y^{-1}}^\vee = {\Cal L}_\Delta$,
cf. (5.23). Therefore, we deduce:
$${\Cal M}_y * {\Cal R} = {\Cal M}_y * 
({\Cal M}_{y^{-1}}^\vee * {\Cal M}_{yw_0}^\vee) =
({\Cal M}_y * {\Cal M}_{y^{-1}}^\vee ) * {\Cal M}_{yw_0}^\vee =
{\Cal L}_\Delta * {\Cal M}_{yw_0}^\vee = {\Cal M}_{yw_0}^\vee\,.$$
It follows that the functor of $\ast$-convolution with ${\Cal R}$
takes modules with ${\Cal M}_y$-flag to modules (not
just complexes in the derived category)
with ${\Cal M}_{yw_0}^\vee$-flag and, moreover, the multiplicities in the two
flags correspond. Since ${\Cal P}_w$ has a ${\Cal M}_y$-flag, we
deduce that ${\Cal Q}_w = {\Cal P}_w * {\Cal R}$ is an actual
module. Furthermore, ${\Cal Q}_w$ has a ${\Cal M}_{yw_0}^\vee$-flag, and we
have: $[{\Cal Q}_w : {\Cal M}_{yw_0}^\vee]= [{\Cal P}_w : {\Cal M}_y]\,.$
This, together with (6.9) and the remark preceding the proof of the theorem,
implies the multiplicity formula of part (iii).

Finally, we observe  that the functor of
$*$-convolution with ${\Cal R}$ is an equivalence 
of derived categories, for it has $*$-convolution with ${\Cal R}_!$
as its inverse.
Therefore, we have: $\Hom_{{\Cal H}{\Cal C}_0}({\Cal Q}_w, {\Cal Q}_y)
\simeq \Hom_{{\Cal H}{\Cal C}_0}({\Cal P}_w, {\Cal P}_y)$.
The last equation of part (iii) is proved similarly
(cf. discussion 
preceding the proof of the theorem).
\qed\enddemo\medskip

\noindent
{\bf {Remarks. }} (i) The multiplicity formula of Theorem 6.9(iii) may be 
viewed as a character formular for tilting modules. An analogous
 formula for tilting modules over an affine Lie algebra is
the main result of [S2]. We observe further
that, in the affine setup there are two affine
flag manifolds: ${\Cal B}_+$ and ${\Cal B}_-$,
with strata of finite co-dimension and strata of finite dimension,
respectively. Category ${\Cal O}$ at a positive level has enough
projectives
and corresponds
 $\widetilde{\Cal D}$-modules on ${\Cal B}_+$, while
category ${\Cal O}$ at a negative level 
corresponds
 $\widetilde{\Cal D}$-modules on ${\Cal B}_-$. The latter category
has {\it no}
projectives, but has well-defined tilting modules instead.
In the affine setup, the "kernel" ${\Cal R}$ lives naturally
on ${\Cal B}_+ \times {\Cal B}_-$.
Our formula: ${\Cal Q}_w
:= {\Cal P}_w * {\Cal R}$ says that the tiltings on ${\Cal B}_-$
are obtained from the projectives on ${\Cal B}_+$
by convolving the latter with ${\Cal R}$, the most natural
way to "transport" projectives from ${\Cal B}_+$ to
${\Cal B}_-$.

(ii) It would be very interesting to prove an affine analogue
of Corollary 6.4 saying that :$\quad$ 
$\Gamma({\Cal B}_+ \times {\Cal B}_-\,,\,{\Cal
R})$
is the {\it semi-infinitely induced} module that played a crucial role
in [S2].

(iii) The same argument as the one used in the proof of
the $\Hom$-equality of Theorem 6.10(iii) also yields:
$$Ext^i_{{\Cal H}{\Cal C}_0}({\Cal Q}_y,{\Cal Q}_w) =
Ext^i_{{\Cal H}{\Cal C}_0}({\Cal P}_y,{\Cal P}_w) = 0\quad,\quad
\forall i>0\,.\quad\qed$$

\enddocument